\documentclass[12pt]{article}
\usepackage{authblk}
\usepackage[spanish,english]{babel}
\usepackage{orcidlink}
\usepackage{amsfonts}
\usepackage{amsmath, hhline, epsfig, latexsym}
\usepackage{amsthm}
\usepackage{amssymb}
\usepackage{mathrsfs}
\usepackage{mathtools}
\mathtoolsset{showonlyrefs}
\usepackage{amsmath}
\numberwithin{equation}{section}
\usepackage{hyperref}
\usepackage{xcolor}
\usepackage[all]{xy}
\usepackage{fancyhdr}
\usepackage{graphicx}
\usepackage{comment}
\usepackage{hyperref}
\usepackage{verbatimbox}
\usepackage{enumerate}
\usepackage[normalem]{ulem}
\usepackage{makeidx} 
\usepackage{subeqnarray}
\usepackage{xcolor}

\newtheorem{theorem}{Theorem}[section]

\def\dis{\displaystyle}

\theoremstyle{definition}
\newtheorem{definition}{Definition}
\newtheorem{remark}{Remark}
\newtheorem{lema}{Lemma}[section]
\newtheorem{proposition}{Proposition}[section]
\usepackage[new]{old-arrows}
\usepackage[nottoc]{tocbibind}
\usepackage{imakeidx}
\usepackage{amsfonts} 
\usepackage{amsmath, hhline, epsfig, latexsym}
\usepackage{amssymb}
\usepackage[all]{xy}
\usepackage{fancyhdr}
\usepackage{graphicx}
\usepackage{hyperref}
\usepackage{enumerate}
\usepackage[normalem]{ulem}
\usepackage{makeidx} 
\usepackage{subeqnarray}
\usepackage{xcolor}
\usepackage[top=1in, left=0.7 in, right=0.6 in,bottom=1in]{geometry}
\usepackage[pagewise]{lineno}
\usepackage[utf8]{inputenc}
\providecommand{\keywords}[1]
{
  \small	
  \textbf{\textit{Keywords---}} #1
}
\providecommand{\subjclass}[1]
{
  \small	
  \textbf{\textit{MSC numbers---}} #1
}

\title{A Stackelberg strategy for the semi-linear heat equation with boundary controls.}

\author{Jose Antonio Villa \orcidlink{0000-0002-9643-3294}\thanks{Published in Boletín de la Sociedad Matemática Mexicana on January 16, 2026: \url{https://doi.org/10.1007/s40590-026-00853-w}}}

\date{}

\begin{document}

\maketitle
\begin{abstract}
    This work continues the results presented in \cite{calsavara2022new} where a hierarchical control problem is solved. Our first goal is the case of both controls acting on the boundary. The follower control solves a null controllability problem, and the leader minimizes a $L^2$ cost functional. Additionally, we solve the case where one control acts on the boundary and the other on an open subset of the domain. More precisely, we prove the existence and uniqueness of a leader-follower couple and give some techniques to solve the problem.

\end{abstract}



\section{Introduction and main results}

This paper addresses a control problem for a semi-linear heat equation in which two controls act on the system with different objectives and with a hierarchical structure. More precisely, the so-called leader imposes their control, while the follower reacts to the leader's actions in order to satisfy its own objective. Classical results on such hierarchical control problems can be found in \cite{araruna2018hierarchic, araruna2020hierarchical}.

In recent work \cite{calsavara2022new}, a hierarchical control strategy was proposed for the semi-linear heat equation in which the classical roles of leader and follower are reversed. Specifically, in \cite{calsavara2022new}, the leader control optimizes a functional, while the follower aims for a null control objective.

Since at least one control acts on the boundary, higher regularity of the solution is required to obtain results in the semi-linear case, introducing additional technical challenges. The main tool to solve this problem is Schauder’s fixed point theorem, where compactness plays a crucial role.

Let  $\Omega\subset \mathbb{R}^n$ be an open set with boundary $\Gamma$ that is a $C^2$ submanifold\footnote{In \cite{LionsMagenes} the boundary of $\Omega$ is an $n-1$ smooth submanifold for simplicity when working with boundary value problems. It is possible to work with more general Lipschitz domains $C^{1,1}$, but here take $C^2$. }. Let $\sigma,\gamma \subset \Gamma$   be relatively open subsets in the relative topology with $\gamma \cap \sigma =\emptyset$. Denote by $Q=\Omega \times (0,T)$ and $\Sigma := \Gamma \times (0,T)$.  Consider an initial datum $y_0\in L^2(\Omega)$  and let $F\in C^2(\mathbb{R})$  globally Lipschtiz function with $F(0)=0$   . Define the initial value problem 

\begin{equation}\label{syst56ghj}
				\left\{\begin{array}[c]{lll}
						\displaystyle y_t-\Delta y +F(y)=0  	& \hbox{in} &Q, \\
                        y = f1_{\gamma}+v1_{\sigma} & \hbox{on} & \Sigma,\\
                        y(0) = y_0 & \hbox{in} & \Omega,
 				\end{array}
                \right.
		\end{equation}
  where $v$ and $f$ belong to suitable Hilbert spaces of functions.   For weight functions blowing up as $t\rightarrow T$, $\varrho \in C^\infty(\bar{Q})$, $\varrho_0 \in C^\infty(\Sigma)$ and $\varrho_*\in C^\infty(0,T)$ to be fixed later, we define

\begin{equation}\label{1.|} 				\mathcal{Y}=\{ y : \varrho y\in L^2(Q)\},\quad\mathcal{F}=\{ f : \varrho_0f\in L^2(\gamma \times (0,T))\},\quad
                \mathcal{V}= \{ v : \varrho_*v\in L^2(\sigma\times (0,T))\}.
 		\end{equation}
  Each set, endowed with the natural  $L^2$ norm gives the Hilbert spaces $(\mathcal{Y},\|\cdot\|_{\mathcal{Y}})$, $(\mathcal{F},\|\cdot\|_{\mathcal{F}})$ and $(\mathcal{V},\|\cdot\|_{\mathcal{V}})$.  

  \begin{definition}
      Define the real function on $\mathcal{V}\times \mathcal{F}$ by
     
           \begin{equation}\label{2..x}
				S(v;f):=\frac{1}{2}\int_Q\varrho^2|y|^2\,dxdt+\frac{1}{2}\int_{\gamma\times(0,T)}\varrho_0^2|f|^2\,d\Sigma.
		\end{equation}
       where $y$ is a solution to \eqref{syst56ghj} that depends on $v$ and $f$ by \eqref{syst56ghj}.
  \end{definition}

 \begin{definition}
    Let  $Q_d:=\Omega_d\times (0,T)$ where   $\Omega_d\subset \Omega$ is an open set of $\mathbb{R}^n$. Take a fixed target $y_d\in L^2(Q_d)$.  On the space $\mathcal{V}\times \mathcal{F}$ define 
      \begin{equation}\label{func76638rtgh}
				P(v;f)=\frac{\alpha}{2}\int_{Q_d}|y-y_d|^2+\frac{1}{2}\int_{\sigma\times(0,T)} \varrho_*^2|v |^2\,d\Sigma.
		  \end{equation} 
     \end{definition}

     \begin{definition}
         We say that $f\in \mathcal{F}$ solves the null controllability problem for $T>0$ if 
         \begin{equation}\label{3.2.5}
             y(T)=0, \hspace{0.3cm} a.e,\hspace{0.2cm} x\in \Omega.
         \end{equation}
     \end{definition}
    Observe that the real function $F$ is required to satisfy $F(0)=0$ to get that the trivial solution $y=0$ solves the null controllability problem \ref{3.2.5}.  Now we present a hierarchical control strategy:\\
 
 \textbf{Stackelberg control strategy:}
 \begin{enumerate}
  \item Given a leader control $v\in \mathcal{V}$  find some follower control $f[v]\in\mathcal{F}$ that solves the null controllability \eqref{3.2.5}.
 
 \item We look for an admissible control $\hat{v}\in \mathcal{V}$  that solves the optimization problem
        
        			\begin{equation}\label{5.6.7}
        							P(\hat{v};f[\hat{v}])=\inf_{v\in \mathcal{V}}P(v;f[v]).
        			\end{equation}

 \end{enumerate}

\begin{remark}
Observe that if for some pair $(v,f[v])\in \mathcal{V}\times \mathcal{F}$, the functional fulfills $S(v;f[v])<\infty$ , then  $y\in \mathcal{Y}$  and blow up at infinity of $\varrho$ when $t\rightarrow T$   implies that $y(T)=0$ a.e on $\Omega$.  If $(v,f[v])\in \mathcal{V}\times \mathcal{F}$ solves the optimization problem 

\begin{equation}\label{4.5.6}
    S(v;f[v])=\inf_{f\in \mathcal{F}} S(v;f),
\end{equation}
then $f[v]$ solves the null controlability problem. Finding a solution $f[v]\in \mathcal{F}$ to problem \eqref{4.5.6} implies the existence of a follower control that solves the null controllability problem \eqref{3.2.5}.
\end{remark}

We are also interested in a variation of the problem \eqref{syst56ghj} where the domain of the controls are changed. Let $\omega \subset \Omega$  such that $\gamma \cap  \omega = \emptyset$. Define the initial value problem

\begin{equation}\label{syst123}
				\left\{\begin{array}[c]{lll}
						\displaystyle y_t-\Delta y +F(y)=v1_{\omega} 	& \hbox{in} &Q, \\
                        y = f1_{\gamma}  & \hbox{on} &\Sigma,\\
                        y(0) = y_0 & \hbox{in} & \Omega.
 				\end{array}
                \right.
		\end{equation}
Consider
	\begin{equation}
               \tilde{ \mathcal{V}}= \{ v : \varrho_0v\in L^2(\omega\times(0,T))\}.
 		\end{equation}
      
Endow this set with the natural weighted $L^2$ norm and define the Banach space $(\tilde{\mathcal{V}},\|\cdot\|_{\tilde{\mathcal{V}}})$. Take $\mathcal{Y}$ and $\mathcal{F}$ as in \eqref{1.|}.

\begin{definition}\label{45.6}
    Let $Q_d:=\Omega_d\times (0,T)$ where $\Omega_d\subset \Omega$ is an open set of $\mathbb{R}^n$ and $y_d\in L^2(Q_d)$.  Define 
      \begin{equation}\label{func76638rtgh1}
				P(v;f)=\frac{\alpha}{2}\int_{Q_d}|y-y_d|^2+\frac{1}{2}\int_{\omega\times(0,T)} \varrho_*^2|v |^2\,d\Sigma.
		  \end{equation} 
     \end{definition}
From the viewpoint of applications, the hypothesis of acting on the boundary arises naturally when some  temperature  is required in some region inside a material body and a heat source is applied on some area of the body borders. Also, it is essential to optimize resources in order to get an objective, that is \textbf{\textit{minimize heat loss}} on an industrial process  (see   \cite{badescu2017optimal}). This idea is formalized by minimizing the cost functional \eqref{func76638rtgh}. With a single heat source, it could not be possible to reach both objectives, therefore, two heat sources \textbf{\textit{ controls}} are required, interacting via the Stackelberg control strategy.

Observe that control by trajectories can be reduced to a null controlability problem. Consider the solution to 
  \begin{equation} \label{a.9-15}
\left\{
\begin{array}[c]{lll} \displaystyle
\bar{y}_{t}-\Delta \bar{y} + F(\bar{y})=0 & \hbox{in} &Q,
\\  \displaystyle
\bar{y} = 0 &\hbox{on}& \Sigma,
\\  \displaystyle
\bar{y}(0) = y_0& \hbox{in}&\Omega.
\end{array}
\right.
   \end{equation}
Change variables $z=y-\bar{y}$, where $y$ is the solution  to \eqref{syst56ghj}. Define the function
\begin{equation}
    H(x)=\int_{0}^1F'(\bar{y}+sx)ds
\end{equation}
with the propriety $G(x)=H(x)x =F(\bar{y}+x)-F(\bar{y})$ . System \eqref{syst56ghj} with the objective 
\begin{equation}
     y(T)=\bar{y}(T) \, \hspace{0.5cm} x\in \Omega,
\end{equation}
can be reformulated as
\begin{equation}\label{mmm}
				\left\{\begin{array}[c]{lll}
						\displaystyle z_t-\Delta z +G(z)=0  	& \hbox{in} &Q, \\
                        z = f1_{\gamma}+v1_{\sigma} & \hbox{on} & \Sigma,\\
                        z(0) = 0, z(T)=0& \hbox{in} & \Omega.
 				\end{array}
                \right.
		\end{equation}
 Define $z_d=y_d-\bar{y}$. The optimization problem for the leader can be formulated with the functional 
\begin{equation}\label{123}
				P(v;f)=\frac{\alpha}{2}\int_{Q_d}|z-z_d|^2\,dxdt+\frac{1}{2}\int_{\sigma\times(0,T)} \varrho_*^2|v |^2\,d\Sigma.
		  \end{equation}

\subsection{Main results and plan of the paper.}
The main results of the paper are presented in two theorems. The first theorem concerns equation \eqref{syst56ghj} where both controls act on the boundary while the second one addresses \eqref{syst123} where the leader is an internal control.

\begin{theorem}\label{1.23}
Let $T>0$ and $y_0\in L^2(\Omega)$ be given. There exists a pair of controls $(\hat{v},\hat{f}[\hat{v}])\in \mathcal{V}\times \mathcal{F}$ such that $\hat{v}$ solves the optimization problem \eqref{5.6.7} and the associated solution $\hat{y}\in \mathcal{Y}$ to \eqref{syst56ghj} satisfies $\hat{y}(T)=0$. Moreover we have the estimate
 \begin{equation}\label{3,.,d1}
 \|\hat{f}[\hat{v}]\|_{\mathcal{F}}+\|\hat{y}\|_{\mathcal{Y}}\leq C\left(\|\hat{v}\|_{\mathcal{V}}+\|y_0\|_{L^2(\Omega)}\right). 
		\end{equation}
\end{theorem}

The first step is  solve the hierarchical control problem for the linear case $F(y)=ay$ with $a\in L^\infty(Q)$ which will be done in Proposition \ref{0.zq}. The strategy is to construct the follower control finding a solution to \eqref{4.5.6} encouraged by Remark 1. \\

In Proposition \ref{5.e.f} we proof the existence of the follower control for the semi-linear case  based on the results of Proposition \ref{0.zq} and Schauder fixed point theorem. Moreover via Dubovitsky- Milyoutin formalism for optimization with constrains, we will find explicit solutions for the follower control for the semi-linear case. The last step is to construct the leader control which is done in  Proposition \ref{OptS.s}.  To verify the existence of the leader that satisfies \eqref{5.6.7} we need Lemma \ref{pp674hg} to assure that the functional $P$ can be minimized in the non linear case.

Now we present the second main result of this paper.

\begin{theorem}\label{1.23T}
Let $T>0$ and $y_0\in L^2(\Omega)$ be given. There exists a pair of controls $(\hat{v},\hat{f}[\hat{v}])\in \tilde{\mathcal{V}}\times \mathcal{F}$ such that $\hat{v}$ solves the optimization problem \eqref{5.6.7} and the associated solution $\hat{y}\in \mathcal{Y}$ to \eqref{syst123} satisfies $\hat{y}(T)=0$. 
\end{theorem}

The proof for this theorem is done following the same strategy of Theorem \ref{1.23}.  The proof of Proposition \ref{0.2zq} and \ref{7.tg}  will not be done with all the explicit computations since proceed in a similar way to Theorem \ref{1.23}, but technical details will be addressed. 

\subsection{Preliminaries}
In this section we recall some  basic results about regularity and compactness of Sobolev spaces. Given $X$ a Banach space, $s\in \mathbb{R}$, $1\le p<\infty$ we define

$$W^{s,p}(0,T;X)=\left\{f\in L^p(0,T;X) \text{ and } \int_0^T\!\!\!\int_0^T\frac{\|f(t)-f(\tau)\|^p_X}{|t-\tau|^{sp+1}} dtd\tau <\infty\right\}.$$
We recall the following compactness result due to Simon \cite{simon1986compact}, (Corollary 9, p. 90).
 \begin{proposition}\label{5.x.}
        	Let  $X,B$ and $Y$ Banach spaces. Assume that $X \subset B$ is compact and $B\subset Y$. For $s_0,s_1$   reals, $\theta\in (0,1)$ and $1\leq r_0\leq \infty,1\leq r_1\leq \infty$, define the numbers $s_{\theta}=(1-\theta)s_0+\theta s_1$, $ \frac{1}{r_\theta}=\frac{\theta}{r_1}+\frac{1-\theta}{r_0}$ and $s_*=s_{\theta}-\frac{1}{r_\theta}$.  Let $A$ be a bounded set in  $W^{s_0,r_0}(0,T;X)\cap W^{s_1,r_1}(0,T;Y)$. If $s_*\leq 0$ then $A$ is relatively  compact in $L^p(0,T;B)$ for $p< - \frac{1}{s_*}$. 
       \end{proposition}
 
 Given an open set $\Omega\subset \mathbb{R}^n$, define the Sobolev space for $r,s$ real numbers 
   		\begin{equation}
   				H^{r,s}((0,T)\times \Omega)=L^2(0,T;H^r(\Omega))\cap H^s(0,T;L^2(\Omega)).
   		\end{equation}
We recall the following result which can be found in \cite{LionsMagenes}, Theorem 5.2, p. 81. This theorem is more illustrative with example 1, Section 15.3: Examples (I) in the same page. 

   \begin{proposition}\label{li..-n}
   		Let $a\in L^\infty(Q)$, $h\in L^2(Q)$, $g\in L^2(\Sigma)$ and $y_0\in L^2(\Omega)$. Then it exists a unique solution $y\in H^{1/2,1/4}(Q)$ that solves the initial value problem
        	\begin{equation}
        			\left\{	\begin{array}[c]{lll}
        						y_t-\Delta y +ay =h & \hbox{in} & Q,\\
                                y=g &\hbox{on} & \Sigma,\\
                                y(0)=y_0 & \hbox{in} & \Omega.    
        			\end{array}\right.
        	\end{equation}
           with the estimate
            \begin{equation}\label{h12}
            				\|y\|_{H^{1/2,1/4}(Q)}\leq C(\|y_0\|_{L^2(\Omega)}+\|h\|_{L^2(Q)}+\|g\|_{L^2(\Sigma)}).
            		\end{equation}
    \end{proposition}   

\begin{proposition} Let $\Omega\subset \mathbb{R}^n$ be a bounded set of class $C^2$. Then it exist a function $\tilde{\eta}_0$ 	
   \[
\displaystyle \tilde\eta_0 \in C^2(\overline{\Omega}), \ \ \tilde \eta_0 \ge 0, \quad\  \nabla \tilde \eta_0\not =0 \ \text{ in }\Omega, \text{ and } 
\frac{\partial\tilde \eta_0}{\partial\nu}\le 0 \ \text{ on } \ \partial \Omega \setminus \gamma .
   \]

\end{proposition}
This classical result  can be found in   Lemma~1.1, p.~4 ~\cite{FI}.  Take $\lambda>0$ and consider $\ell \in C^2([0,T])$ that satisfies
$
\ell(t) \geq T^2/4 \ \text{ in } \ [0,T/2] \ \text{ and } \ \ell(t) = t(T-t) \ \text{ in } \ [T/2,T].
$ We introduce the following    functions  
   
   \begin{equation}\label{weightk}
\begin{array}{c} \displaystyle
{\displaystyle \tilde{\sigma}(x,t):=\frac{e^{4\lambda \|\tilde{\eta}^{0}\|_{\infty}}-e^{\lambda(2\|\tilde{\eta}^{0}\|_{\infty}+\tilde{\eta}^{0}(x))}}{\ell(t)}},\ \ \ {\displaystyle\tilde{\xi}(x,t):=\frac{e^{\lambda (2\|\tilde{\eta}^{0}\|_{\infty}+\tilde{\eta}^{0}(x))}}{\ell(t)}}.
\end{array}
   \end{equation}
   
 We consider now a parameter $s>0$ and define the weights~$\varrho=e^{s\tilde{\sigma}}$, $\varrho_{0} = (s\tilde\xi)^{-3/2} \lambda^{-2}\varrho $, $\varrho_{1} =(s\tilde\xi)^{-1/2}\lambda^{-1}  \varrho $,  $\varrho_{2} =(s \tilde\xi)^{1/2} \,\varrho$. 
 The parameters $\lambda$ and $s$ are going to be chosen large enough.  Define the vector space $\mathcal{P}_0=\{q\in C^2(\overline{Q}):q|_{\Sigma}=0\}$.
\begin{proposition}[Carleman Inequality]\label{carlm3}
   There exist positive constants $\lambda_0$, $s_0$ and~$C_1$, only depending on $\Omega$, $\gamma$ and~$T$, such that, if we take $\lambda=\lambda_0$ and~$s= s_0$, any~$p \in \mathcal{P}_0$ satisfies
   \begin{equation}\label{mnb}
\int_{Q} \left[ \varrho_{2}^{-2}(|p_t|^2 + |\Delta p|^2) + \varrho_{1}^{-2} |\nabla p|^2 + \varrho_{0}^{-2} |p|^2 \right]
\leq C_1 \left( \int_{Q}\varrho^{-2}|p_t+\Delta p|^2+\int_{\gamma\times (0,T)}\varrho_0^{-2}|\partial_{\eta}p|^2\,d\Sigma  \right).
   \end{equation}
Furthermore, $\lambda_0$ and~$s_0$ can be found arbitrarily large.
\end{proposition}

For $a\in L^\infty (Q)$ consider the operator $L_a=\partial_t-\Delta  +a$ and its adjoint   $L_a^*=-\partial_t-\Delta +a$.
Define the bilinear form $B:\mathcal{P}_0\times \mathcal{P}_0\longrightarrow \mathbb{R}$  by 
 
	\begin{equation}\label{8.gf}
			B(a;p,q)=\int_{Q}\varrho^{-2}L_a^*(p)L_a^*(q)\,dxdt+\int_{\gamma\times (0,T)}\varrho_0^{-2}\partial_{\eta}p\partial_{\eta}q\,d\Sigma.
	\end{equation}
 The above inequality implies that  $B$ is a coercive bilinear form. Define the semi-norm  $\|p\|_{\mathcal{P}_0}:=B(0;p,p)^{1/2}$.  Carleman inequality \eqref{mnb} implies $\|p\|_{\mathcal{P}_0}=0$ then $\varrho_0^{-2}|p|^2=0$ a.e. so  $\|\cdot\|_{\mathcal{P}_0}$ is a norm in $\mathcal{P}_0$.  Define $\mathcal{P}$ the completion of $\mathcal{P}_0$ with the norm $\|\cdot\|_{\mathcal{P}_0}$.   We can find positive constants $K_0$ and~$K_1$, only depending on $\Omega$, $\gamma$, $T$ and~$\|a\|_{L^\infty(Q)}$, such that

   \begin{equation}\label{6,4.}
K_0 \, B(0;p,p) \leq B(a;p,p) \leq K_1 \, B(0;p,p), \quad \forall p \in \mathcal{P}.
  \end{equation}
 \begin{lema}\label{7.56t}
The map $\mathcal{P}\longrightarrow L^2(\Omega)$, $p\mapsto p(0)$ is continuous.
 \end{lema}
 
 \begin{proof}
  Take $\delta >0$. From Carleman inequality we see that $p1_{[0,T-\delta]}\in L^2(0,T;H^2(\Omega))$ and the derivative $p_{t}1_{[0,T-\delta]}\in L^2(0,T/2;L^2(\Omega))$. By Theorem 3.2 pg. 21 \cite{LionsMagenes},  $p(0)\in H^1_0(\Omega)$ and this trace map is continuous.
  \end{proof} 
 
\noindent   Given  $p\in \mathcal{P}$, 
  
  		\begin{equation}\label{m..-4}
  				\int_{Q}\varrho_*^{-2}|\Delta p|^2 +\varrho_*^{-2}|\nabla p|^2+\varrho_*^{-2}|p|^2\,dxdt \leq \int_{Q}\varrho_2^{-2}|\Delta p|^2 +\varrho_1^{-2}|\nabla p|^2+\varrho_0^{-2}|p|^2\,dxdt < \infty.
  		\end{equation}
      Then,  $\Delta(\varrho_*^{-1}p),\nabla(\varrho_*^{-1}p), \varrho_*^{-1}p\in L^2(Q)$ and   $\varrho_*^{-1}p\in L^2(0,T;H^2(\Omega))$.

      		\begin{lema}\label{4--6}
      				Let 
                      $\sigma\subset \Gamma$ and open set. Then, for any $p\in \mathcal{P}$ it exists a constant $C>0$ such the next inequality holds
                    \begin{equation}\label{3,bn}
                    		\int_{\sigma\times(0,T)}\varrho_*^{-2}|\partial_\eta p|^2\,d\Sigma \leq C \|p\|^2_{\mathcal{P}}.
                            \end{equation}
      		\end{lema}
            
          \begin{proof}
         Given $p\in \mathcal{P}$, we have from Carleman Inequality \eqref{m..-4}, that  $\varrho^{-1}_*p\in L^2(0,T;H^2(\Omega))$. From Theorem 1.1,  pg. 9 \cite{lions2012non}, the operator $\partial_\eta:L^2(0,T;H^2(\Omega))\longrightarrow L^2(0,T;H^{1/2}(\Gamma))$ is continuous. Then,
          
          		\begin{equation}
          				\int_{\Sigma}|\partial_\eta (\varrho^{-1}_*p)|^2\,d\Sigma \leq \tilde{C}\int_{Q}\varrho_*^{-2}|\Delta p|^2 +\varrho_*^{-2}|\nabla p|^2+\varrho_*^{-2}|p|^2\,dxdt.
          		\end{equation}
          
Combining \eqref{m..-4} and Carleman inequality  
          		\begin{equation}
          				\int_\Sigma \varrho_*^{-2}|\partial_\eta p|^2\,d\Sigma \leq C \|p\|^2_{\mathcal{P}}.
          		\end{equation}
                 Finally since $\sigma\subset \Gamma$, then \eqref{3,bn} holds.
          \end{proof}

  \begin{proposition}[Dubovitski-Milyoutin][\cite{alexeev2017commande}]\label{dv.yu}
      Let $H$ a Hilbert space, $F$  a real function on $H$ and $M\in C^1(H,H)$. Let $\hat{h}\in H$ be a solution to the optimization problem
      \begin{equation}
          F(\hat{h})=\inf_{h\in H} F(h), \hspace{0.5cm}M(\hat{h})=0.
      \end{equation}
      Then it exist a $\lambda >0$ and $\zeta\in \text{ker}(M'(\hat{h}))^\bot$     not both zero such that 
      \begin{equation}\label{1.3gaw}
         \lambda F'(\hat{h})+\zeta=0.
      \end{equation}
Observe that if $M$ is a compact operator then $M'$ is compact and then $$\ker M'(\hat{h})^\bot=\text{Rank} (M'(\hat{h}))^*.$$    
  \end{proposition}

\section{Boundary  controls. }
In this section we  prove of Theorem \ref{1.23}. First, we solve the null controlability problem for the linear case and then, in subsection 2.1, we solve the semi-linear problem. For sake of  completeness we present the proof of the linear case as done in section 5.3 in \cite{calsavara2022new}.

\begin{proposition}\label{0.zq}
 Let $v\in \mathcal{V}$,  $T>0$, $y_0\in L^2(\Omega)$ and  $a\in L^\infty(Q)$ given. It exists a follower control $f[v]\in \mathcal{F}$ such that the solution $y\in \mathcal{Y}$ of the initial value problem
   		\begin{equation}\label{2.5xc}
					\left\{\begin{array}[c]{lll}
						\displaystyle y_t-\Delta y +ay=0 	& \hbox{in} &Q, \\
                        y = f[v]1_{\gamma}+v1_{\sigma}  & \hbox{on} &\Sigma,\\
                        y(0) = y_0 & \hbox{in} & \Omega,
 				\end{array}\right.
		\end{equation}
satisfies $y(T)=0$ and
 
 \begin{equation}\label{8.xy}
 f[v]=\varrho_0^{-2}\partial_\eta p1_{\gamma}; \hspace{0.5cm} y=\varrho^{-2}L_a^*(p),
 \end{equation}
 where $p\in \mathcal{P}$  solves
 	\begin{equation}\label{v.z}
				B(a;p,q) = \int_{\sigma\times(0,T)}v \partial_{\eta}q\, d\Sigma+\langle y_0,q(0) \rangle_{L^2(\Omega)},\hspace{0.2cm} \forall q\in \mathcal{P}. 
		\end{equation}
 Moreover 
 \begin{equation}\label{3,.,d}
 \|f[v]\|_{\mathcal{F}}+\|y\|_{\mathcal{Y}}\leq C\left(\|v\|_{\mathcal{V}}+\|y_0\|_{L^2(\Omega)}\right). 
		\end{equation}
\end{proposition}

\begin{proof}
\textbf{Step 1}.  The functional $S(v;\,\cdot ):\mathcal{F}\longrightarrow \mathbb{R}$ is a coercive, convex, and weakly lower semicontinous functional. Then, it has a unique minimizer $f[v]\in \mathcal{F}$ that solves problem \eqref{4.5.6}. Define the functions $p$ and $z$  solutions to 
        \begin{equation}\label{..--c}
 					\left\{\begin{array}[c]{lll}
 						\displaystyle -p_t-\Delta p +ap =\varrho^2 y &\hbox{in} & Q,\\
                        p=0 &\hbox{on} & \Sigma,\\
                        p(T)=0 &\hbox{in} & \Omega,
 				\end{array}\right.
 		\end{equation}
        and for $h\in \mathcal{F}$,     
        		\begin{equation}\label{.bn-c}
 					\left\{\begin{array}[c]{lll}
 						\displaystyle z_t-\Delta z +az =0 &\hbox{in} & Q,\\
                        z=h1_{\gamma} &\hbox{on} & \Sigma,\\
                        z(0)=0 &\hbox{in} & \Omega.
 				\end{array}\right.
 		\end{equation}
Let now compute the Gateaux derivative of $S(\cdot;v)$ at $f[v]$:
 				\begin{equation}\label{7.,x}
   				 \lim_{\epsilon \to 0}\frac{1}{\epsilon}\left[S(v; f[v]+\epsilon h)-S(v;f[v]))\right] = \int_{Q}\varrho^2yz\,dxdt+\int_{\gamma\times (0,T)}\varrho^2_0f[v]h d\Sigma=0.
   		\end{equation}
   From \eqref{7.,x}, multiplying \eqref{..--c} by $z$ and integrating in $Q$,  we get        

            \begin{equation}
            		\int_{\gamma\times(0,T)}(\varrho^2_0f[v]-\partial_\eta p)hd\Sigma=0, \hspace{0.3cm} \forall h\in \mathcal{F}.
            \end{equation}
 This implies the characterization \eqref{8.xy}.

\textbf{Step 2.} Observe that \eqref{8.xy} allows to write  problem \eqref{2.5xc} in the form
	
    	\begin{equation}\label{2.,c}
				\left\{	\begin{array}[c]{lll}
						L_a(\varrho^{-2}L_a^*(p))=0 	& \hbox{in} &Q, \\
                        \varrho^{-2}L_a^*(p) = \varrho_0^{-2}\partial_\eta p1_{\gamma}+v1_{\sigma}  & \hbox{on} &\Sigma,\\
                        \varrho^{-2}L_a^*(p)(0) = y_0; \varrho^{-2}L_a^*(p)(T)=0 & \hbox{in} & \Omega.
 				\end{array}\right.
		\end{equation}
Multiply equation \eqref{2.,c} by $q\in \mathcal{P}$, integrate by parts to get 
		\begin{equation}\label{v.tz}
				B(a;p,q)= \int_{\sigma\times(0,T)}v\partial_{\eta}q\,d\Sigma+\langle y_0,q(0) \rangle_{L^2(\Omega)}. 
		\end{equation}
   From  Theorem \ref{carlm3},we see that the left hand side of equation \eqref{v.tz} is coercive. It remains to verify the continuity of 
\begin{equation}\label{34.k}
    \ell(q) =\int_{\sigma\times(0,T)}v\partial_{\eta}q\,d\Sigma+\langle y_0,q(0) \rangle_{L^2(\Omega)}.
\end{equation}

By Lemma \ref{7.56t} and Lemma \eqref{4--6}we obtain
    		\begin{equation}\label{continuidad}
    				|\ell(q)|\leq\left(\|v\|_{\mathcal{V}}+\|y_0\|_{L^2(\Omega)}\right)B(0;q,q)^{1/2}.
    		\end{equation}
So $\ell$ is  continuous  and  by Lax-Milgram  Theorem  equation \eqref{v.tz}  has a unique  solution $p\in \mathcal{P}$.  It is not difficult to see   that \eqref{v.tz} is equivalent to \eqref{2.,c}

  \textbf{ Step 3.} This final step is to prove inequality \eqref{3,.,d}. In equation \eqref{v.tz} make $q=p$ to get
    	\begin{equation}\label{9.21}
    			B(a;p,p)= \int_{\sigma\times (0,T)}v\partial_\eta p\,d\Sigma+\langle y_0,p(0) \rangle_{L^2(\Omega)}. 
    	\end{equation}
 Then, it is straightforward to get    
     \begin{equation}\label{5-q}
     B(a;p,p)\leq \left(\int_{\sigma\times(0,T)}\varrho_*^{2}|v|^2\right)^{1/2}\left(\int_{\sigma\times(0,T)}\varrho_*^{-2}|\partial_\eta p|^2\right)^{1/2}+\|y_0\|_{L^2(\Omega)}\|p(0)\|_{L^2(\Omega)}.
     \end{equation}
 Inequality  \eqref{5-q} combined with \eqref{continuidad} implies
    
    	\begin{equation}
 B(a;p,p)\leq\left(\|v\|_{\mathcal{V}}+\|y_0\|_{L^2(\Omega)}\right)B(0,p,p)^{1/2}.
   	\end{equation}
 From inequality   \eqref{6,4.}, for $p\not =0$, it exists a constant $C>0$    such that 
   	\begin{equation}
   			B(a;p,p)^{1/2}\leq C\left(\|v\|_{\mathcal{V}}+\|y_0\|_{L^2(\Omega)}\right),
   	\end{equation}
    and then by construction of $ f[v]$ and $ y$; and Young inequality, we get 
    		\begin{equation}		\| f[v]\|_{\mathcal{F}}+\| y\|_{\mathcal{Y}}\leq C\left(\|v\|_{\mathcal{V}}+\|y_0\|_{L^2(\Omega)}\right).
    		\end{equation}
\end{proof}

  \subsection{The semi-linear case.}
  The null controllability problem for the linear case $F(y)=ay$ with $a\in L^\infty(Q)$ solved in proposition \eqref{0.zq} will allow to solve the corresponding problem for the non linear case by a fixed point argument.  Define 
  
  	\begin{equation}\label{lienali644gh}
    		F_0(s)=
  			\left\{\begin{array}[c]{lll}
  					\displaystyle \frac{F(s)}{s} & s\neq 0, \\
               		 0 & s=0.	
  			\end{array}
            \right.
  	\end{equation}
Given a function $z\in L^2(Q)$, we have that $F_0(z)\in L^\infty(Q)$. Define the linearization of (\ref{syst56ghj})  by
    		\begin{equation}\label{s6.}
				\left\{\begin{array}[c]{lll}
						\displaystyle y_t-\Delta y +F_0(z)y=0 	& \hbox{in} &Q, \\
                        y = f1_{\gamma}+v1_{\sigma}  & \hbox{on} &\Sigma,\\
                        y(0) = y_0 & \hbox{in} & \Omega.
 				\end{array}
                \right.
		\end{equation}
Since $F_0(z)\in L^\infty(Q)$, apply the results of Proposition \ref{0.zq} to system \eqref{s6.} as it is done in the next proposition.

\begin{proposition}\label{5.e.f}
 Let a leader control $v\in \mathcal{V}$ ,  $T>0$  and $y_0\in L^2(\Omega).$ Then, it exists a follower control $f[v]\in \mathcal{F}$   where $y\in \mathcal{Y}$ solves the initial value problem
 							\begin{equation}
 										\left\{\begin{array}[c]{lll}
						\displaystyle y_t-\Delta y +F(y)=0 	& \hbox{in} &Q, \\
                        y = f[v]1_{\gamma}+v1_{\sigma}  & \hbox{on} &\Sigma,\\
                        y(0) = y_0,  y(T)=0& \hbox{in} & \Omega,
 									\end{array}\right.
 							\end{equation}
with
  		\begin{equation}\label{e5.}
   				f[v]=\varrho_0^{-2} \partial_\eta p1_{\gamma}, \hspace{1cm}  y =\varrho^{-2}L^*_{F'(y)}(p),
   		\end{equation}
    where  the function $p\in \mathcal{P}$ solves the equation
    
    		\begin{equation}\label{0--j}
				\int_Q\varrho^{-2}L^*_{F'(y)}(p)L^*_{F_0(y)}(q)+\int_{\gamma \times(0,T)}\varrho_0^{-2}\partial_{\eta} p \partial_{\eta} q\, d\Sigma\\
                =\int_{\sigma \times (0,T)} v\partial_\eta q\,d\Sigma+\int_\Omega y_0(x)q(0)dx,  \hspace{0.3cm} \forall q\in \mathcal{P}.
		\end{equation}
 Also, we get the estimate
 \begin{equation}\label{6-2-}
 \|f[v]\|_{\mathcal{F}}+\|y\|_{\mathcal{Y}}\leq C\left(\|v\|_{\mathcal{V}}+\|y_0\|_{L^2(\Omega)}\right).
 \end{equation}
 
  \end{proposition}
  
  \begin{proof}
\textbf{Step 1.}   Let $z\in L^2(Q)$ be given. Consider  the functions $f_z=f[v]_z$ and $y_z$ given by Proposition \ref{0.zq}. By \eqref{3,.,d}, the follower control is uniformly bounded   as \[\|f_z\|_{\mathcal{F}}\leq C\left(\|v\|_{\mathcal{V}}+\|y_0\|_{L^2(\Omega)}\right), \hspace{0.3cm} \forall z\in L^2(Q).\] Moreover, Proposition \ref{li..-n} implies that, for every $z\in L^2(Q)$, the corresponding solution $y_z$ to the problem
			\begin{equation}\label{90--,}
						\left\{\begin{array}[c]{lll}
							\displaystyle y_{z,t}-\Delta y_z +F_0(z)y_z=0 	& \hbox{in} &Q, \\
                        y_z = f_z[v]1_{\gamma}+v1_{\sigma}  & \hbox{on} &\Sigma,\\
                        y_z(0) = y_0; y_z(T)=0 & \hbox{in} & \Omega,
					\end{array}\right.
			\end{equation}
  fulfills $y_z\in H^{1/2,1/4}(Q)$. Then, from the fact that $\|v\|_{L^2(\Sigma)}\leq \|v\|_{\mathcal{V}}$ and inequality \eqref{h12}
		\begin{equation}
				\|y_z\|_{H^{1/2,1/4}(Q)}\leq C\left(\|y_0\|_{L^2(\Omega)}+\|v\|_{\mathcal{V}}+\|f_z\|_{\mathcal{F}} \right ).
		\end{equation}
Since the embedding $H^{1/2}(\Omega)\hookrightarrow L^2(\Omega)$ is compact, we can apply Proposition \ref{5.x.} with $\theta =1/2$, $s_0=1/2$, $s_1=1/4$ and  $s*=-1/4$ getting that the embedding $H^{1/2,1/4}(Q)\hookrightarrow L^2(Q)$ is compact.  Define the map 
\begin{equation}
    \Lambda:L^2(Q)\longrightarrow L^2(Q), \hspace{0.2cm}z \mapsto y_z
\end{equation}
where $y_z$ solves \eqref{90--,}.  By the previous conclusions   $\Lambda (L^2(Q))\subset H^{1/2,1/4}(Q)$  is bounded and then $\Lambda (L^2(Q)) \subset L^2(Q)$ is compact. By Schauder's fixed point Theorem, it exists a $z:=\tilde{y}$ that solves \eqref{90--,} and since $F_0(\tilde{y})\tilde{y}=F(\tilde{y})$ then
 
 			\begin{equation}\label{64-}
 						\left\{\begin{array}[c]{lll}
 							\displaystyle \tilde{y}_t-\Delta \tilde{y} +F(\tilde{y})=0 	& \hbox{in} &Q, \\
                        \tilde{y} = f_{\tilde{y}}[v]1_{\gamma}+v1_{\sigma}  & \hbox{on} &\Sigma,\\
                       \tilde{y}(0) = y_0; \tilde{y}(T)=0 & \hbox{in} & \Omega.
 					\end{array}\right.
 			\end{equation}
            
Denote by $\tilde{f}[v]=f_{\tilde{y}}[v]$. By construction $\tilde{y}\in \mathcal{Y}$ and $\tilde{f}[v]\in \mathcal{F}$, then it is possible to see that 

			\begin{equation}\label{3.jk}
					S(v; \tilde{f}[v])=\frac{1}{2}\int_{Q}\varrho^2|\tilde{y}|^2\,dxdt+\frac{1}{2}\int_{\gamma\times (0,T)}\varrho_0^2|\tilde{f}[v]|^2\,d\Sigma < \infty.
			\end{equation}

\textbf{Step 2.} We obtain that the set of $f\in \mathcal{F}$  where the function $S(v;f)< \infty$  and $y$ solves the semi-linear problem is non empty. Since \eqref{3.jk} holds, the function $S(v;\,\cdot\,)$   has a positive  lower bound. Then it exists  $\{f_n\}\subset \mathcal{F}$  a uniformly bounded  minimizing sequence for $S(v; \cdot)$.   For a sub-sequence $f_{n_k}$, the associated states $y_n$ converges strongly to $y$ in $L^2(Q)$ and $f_n$ converges weakly to some $\hat{f}$ in $\mathcal{F}$. From the weakly lower semi-continuity of the functional
   \[
S(v; \hat{f})\leq \liminf_{n\rightarrow \infty}S(v;f_n)=\inf_{f\in \mathcal{F}}S(v;f),
   \]
then $\hat{f}\in \mathcal{F}$ is a  solution to problem \eqref{4.5.6}.

\textbf{Step 3.}  Define $\bar{y}\in C(0,T;H^1_0(\Omega))\cap C^1(0,T;L^2(\Omega))$  as  the solution to the system
		\begin{equation}
					\left\{\begin{array}[c]{lll}
						\bar{y}_t-\Delta \bar{y}=0 &\hbox{in} & Q,\\
                        \bar{y}=0 & \hbox{on} & \Sigma,\\
                        \bar{y}(0)=y_0 &\hbox{in} & \Omega.
				\end{array}\right.
		\end{equation}
Define the map $H_0:L^2(Q)\longrightarrow L^2(Q)$ as $H_0(q)=z$ where $z$ is the solution to the problem
			\begin{equation}\label{.fg3}
						\left\{	\begin{array}[c]{lll}
								z_t-\Delta z=q &\hbox{in} & Q,\\
                        z=0 & \hbox{on} &  \Sigma, \\
                        z(0)=0 &\hbox{in} & \Omega.
						\end{array}\right.
			\end{equation}
Observe that $H^*_0:L^2(Q)\longrightarrow L^2(Q)$ is given by $H^*_0(\psi)=\varphi$ where $\varphi$ solves the equation

				\begin{equation}\label{t--g}
							\left\{\begin{array}[c]{lll}
								-\varphi_t-\Delta \varphi =\psi &\hbox{in} & Q,\\
                                \varphi =0  &\hbox{on} &\Sigma,\\
                                \varphi(T)=0  &\hbox{in} & \Omega.
						\end{array}\right.
				\end{equation}
            For each $q\in L^2(Q)$, it holds that  $H_0(q)\in L^2(0,T;H^1_0(\Omega))\cap H^1(0,T;L^2(\Omega))$ . The  map $H_0$ is then a compact operator.   By Proposition \ref{li..-n}, it is possible to define the boundary operator $G:L^2(\Sigma)\to H^{1/2,1/4}(Q)\subset L^2(Q)$ given by $G(\beta)=\eta$ , where $\eta$ is a solution to the boundary value problem
		\begin{equation}
				\left\{	\begin{array}[c]{lll}
						\eta_t-\Delta \eta=0 &\hbox{in} & Q,\\
                        \eta=\beta & \hbox{on} &  \Sigma, \\
                        \eta(0)=0 &\hbox{in} & \Omega.	
				\end{array}\right.
		\end{equation}
  By Proposition \eqref{5.x.}, the operator $G:L^2(\Sigma)\longrightarrow L^2(\Omega)$ is compact. For a fixed $v\in \mathcal{V}$,  define the map $M:\mathcal{Y}\times \mathcal{F}\longrightarrow  L^2(Q)$ given by 

	\begin{equation}\label{mqas}
			M(y,f)=y+H_0(F(y))-G(v1_\sigma+f1_{\gamma})-\bar{y}.
	\end{equation}
    
 It is straightforward to verify that if the condition $M(y,f)=0$ holds, then the pair $(y,f)\in\mathcal{Y}\times \mathcal{F}$ solves the equation 
 		\begin{equation}\label{eq6573f}
 					\left\{\begin{array}[c]{lll}
 						y_t-\Delta y + F(y)=0 &\hbox{in} & Q,\\
                        y = f1_{\gamma} +v1_{\sigma} &\hbox{on} & \Sigma,\\
                        y(0)=y_0 &\hbox{in} &\Omega.
 				\end{array}\right.
 		\end{equation}
        The operators $H_0$ and $G$ are of class $C^1$ so it is the operator $M$. Given any directions $(z,g) \in\mathcal{Y}\times \mathcal{F}$ the derivative of $M$ is the operator $M':\mathcal{Y}\times \mathcal{F}\longrightarrow L^2(Q)$ given by
  		\begin{equation}
  				M'(y,f)(z,g)=z+H_0\left(F'(y)z\right)-G(g1_{\gamma}).
  		\end{equation}
  The operator $M$ is compact since $H_0$ and $G$ are compact operators. Then  $M'(y,f)$ is compact, so is possible to conclude that $\text{ker} M'(y,f)^\bot=\text{Rank} M'(y,f)^*$. It is necessary to compute the dual operator  $M'(y,f)^*:L^2(Q)\longrightarrow \mathcal{Y}\times \mathcal{F}$ given by 
	\begin{equation}\label{s.gh}
			M'(y,f)^*(\psi)=\left(\psi+F'(y)H_0^*(\psi), - G^*(\psi)1_\gamma\right).
	\end{equation}
 The idea is to formulate as an optimization problem with constrains in the form
 
 		\begin{equation}\label{p54t}
 				\left\{
                \begin{array}[c]{lll}
 							\displaystyle\inf \frac{1}{2}\int_{Q}\varrho^2|y|^2 +\frac{1}{2}\int_{\gamma\times (0,T)}\varrho_0^2|f|^2,\\
                            M(y,f)=0,\\
                            (y,f)\in \mathcal{Y}\times \mathcal{F}.
 				\end{array}
                \right.
 		\end{equation}

 Let $(y,f)\in \mathcal{Y}\times \mathcal{F}$  a solution to \eqref{p54t} for a fixed leader $v$. Consider $\lambda >0$ and $(h,k)\in \text{Rank}(M'(y,f)^*)$ given by the Dubovitsky-Milyoutin formalism in Proposition \ref{dv.yu}.  For $\psi\in L^2(Q)$ such that $(h,k)=M'(y,f)^*(\psi)$ observe that, after normalizing $\lambda=1$, \eqref{1.3gaw} takes the form
 	\begin{equation}\label{ehgjy684}
 					(\varrho^2y,\varrho^2_0f[v])+\left(\psi+F'(y)H_0^*(\psi),-G^*(\psi)1_\gamma\right)=0.
 	\end{equation}
  Is possible to deduce that 

    		\begin{equation}\label{|.d}
    				\begin{array}[c]{lll}
    					\dis	y&=\varrho^{-2}\left(\psi+F'(y)H_0^*(\psi)\right),\\
                            \dis f&=-\varrho_0^{-2}G^*(\psi)1_\gamma.
    				\end{array}
    		\end{equation}
Define  $p:=H_0^*(\psi)$ so by definition $L_0^*(p)=L_0^*(H_0^*(\psi))=\psi$ and because $(y,f)\in \mathcal{Y}\times \mathcal{F}$  solves problem \ref{p54t}  then $p\in \mathcal{P} $.   In equation \eqref{|.d} it is possible to write $f[v]=-\varrho_0^{-2}G^*(L_0^*(p))$ and it  is necessary to compute explicitly  $G^*(L_0^*(p))$.  It follows from the next lemma which proof is given in the appendix:

\begin{lema}\label{lemmanm}
    Let $p\in L^2(0,T;H^2(\Omega))$. Then $G^*(L_0^*(p))=-\partial_\eta p$.
\end{lema}

So we have $f[v]=\varrho_0^{-2}\partial_\eta p$. The first equation in equation \eqref{|.d} becomes
        		\begin{equation}
        				\begin{array}[c]{lll}
        						y& = \varrho^{-2}(\psi+F'(y)H_0^*(\psi))\\
                                  &=\varrho^{-2}(L_0(p)+F'(y)p)\\
                                &=\varrho^{-2}L^*_{F'(y)}(p).
        				\end{array}
        		\end{equation}
Then \eqref{e5.} holds. Again, since $(y,f[v])\in \mathcal{Y}\times \mathcal{F}$   fulfills  the restrictions of  Problem  \ref{p54t} then with the above equation it is possible to write  the restriction $M(y,f[v])=0$   in the form
                \begin{equation}\label{p7.64}
                				\left\{\begin{array}[c]{lll}
                				\displaystyle L_{F_0(y)}(\varrho^{-2}L^*_{F'(y)}(p))=0 &\hbox{in} & Q,\\
                                \varrho^{-2}L^*_{F'(y)}(p)=\varrho_0^{-2}\partial_\eta p1_{\gamma}+v1_{\sigma}&\hbox{on} & \Sigma,\\
                                \varrho^{-2}L^*_{F'(y)}(p)(0)=y_0 &\hbox{in} & \Omega.
                		\end{array}\right.
                \end{equation}
\textbf{Step 4}. Take $p'\in \mathcal{P}$ and multiply by it in equation \ref{p7.64}. Using integration by parts and boundary conditions, the integral form for this problem is
	\begin{equation}\label{inhgy65}
        	 \displaystyle \int_Q\varrho^{-2}L^*_{F'(y)}(p)L_{F_0(y)}^*(p')-\int_\Omega y_0(x)p'(0,x)
               -\int_{\sigma \times (0,T)} v(\partial_{\eta}p')\,d\Sigma+\int_{\gamma\times(0,T)}\varrho_0^{-2}\partial_{\eta} p\,\partial_{\eta} p'\,d\Sigma=0.               
		\end{equation}
then from this equations it is possible to get the integral equation taking $p'=p$ 
		\begin{equation}\label{biy576}
				\int_Q\varrho^{-2}L^*_{F'(y)}(p)L^*_{F_0(y)}(p)+\int_{\gamma \times(0,T)}\varrho_0^{-2}|\partial_{\eta} p|^2\, d\Sigma\\
                =\int_{\sigma \times (0,T)} v(\partial_{\eta}p)\,d\Sigma+\int_\Omega y_0p(0)dx . 
		\end{equation}
Now it is necessary to make some estimates from \eqref{biy576}. For this, apply  H\"older and Young inequalities. First it is important to note that the weight $\varrho$ is bounded in the interval $[0,T/2]$ and by Lemma \ref{7.56t}  is possible to see that $\|p(0)\|_{L^2(\Omega)}\leq \max_{t\in [0,T/2]}\|p(t)\|_{L^2(\Omega)}$. Denote $M=\sup_{y\in \mathbb{R}}|F'(y)|$. From Theorem \ref{carlm3} fixing $\lambda>\lambda_0 $ and $s>s_0$ the following inequality holds : 
\begin{equation}
s^{3/2}>M\sqrt{2}\lambda^{-2}\sup_{Q}\xi^{-3/2}.
\end{equation}

Then, define $S:=\sup_Q\frac{\varrho_0}{\varrho}< 1/(M\sqrt{2})$ . It is not difficult to see that

\begin{equation}
\int_{Q}\varrho^{-2}|p|^2dxdt\leq S^2\int_Q\varrho_{0}^{-2}|p|^2\,dxdt.
\end{equation}

Also there exist a positive number $\beta$ such that $\frac{1}{M}\frac{S^2M^2}{1-M^2S^2}<\beta<\frac{1}{M}$. Then 
     
		\begin{equation}\label{ine.75648}
				\begin{array}[c]{lll}
				B(0;p,p)&= \displaystyle \int_\Omega y_0(x)p(0,x)+\int_{\sigma\times(0,T)}v\partial_{\eta}p\,d\Sigma\\
                & \displaystyle- \int_Q \varrho^{-2}\left(F_0(y)L^*_0(p)p+F'(y)L_0^*(p)p+F'(y)F_0(y)|p|^2\right)\,dxdt \\
                & \displaystyle \leq \left( \int_{\sigma\times(0,T)}\varrho_*^{2}|v|^2\,d\Sigma \right)^{\frac{1}{2}}\left( \int_{\sigma\times(0,T)}\varrho_*^{-2}|\partial_\eta p|^2\,d\Sigma \right)^{\frac{1}{2}}+\left( \int_{\Omega}|y_0|^2 \right)^{\frac{1}{2}}\left( \int_{\Omega}|p(0)|^2 \right)^{\frac{1}{2}}\\
                & \displaystyle - \int_Q \varrho^{-2}\left(F_0(y)L^*_0(p)p+F'(y)L_0^*(p)p+F'(y)F_0(y)|p|^2\right)\,dxdt \\
                &  \displaystyle \leq \|y_0\|_{L^2(\Omega)}\max_{t\in[0,T/2]}\|p(t)\|_{L^2(\Omega)} +\|v\|_{\mathcal{V}}\|p\|_{\mathcal{P}}+M^2\int_{Q}\varrho^{-2}|p|^2 \,dxdt\\
                & \displaystyle + 2M\int_Q\varrho^{-2}|L_0^*(p)||p|\,dxdt.  \\               
          \end{array}
          \end{equation}
           Then by Young inequality  with parameter $\beta$ it is possible to bound
           
           					\begin{equation}
           							 2M\int_Q\varrho^{-2}|L_0^*(p)||p|\,dxdt \leq M\beta\int_Q\varrho^{-2}|L^*_0(p)|^2\,dxdt+\frac{M}{\beta}\int_{Q}\varrho^{-2}|p|^2\,dxdt,
           					\end{equation}
Since $\varrho_0< S\varrho$ then we get

			\begin{equation}
						\int_Q\varrho^{-2}|p|^2\,dxdt\leq S^2\int_Q\varrho_0^{-2}|p|^2\,dxdt.
			\end{equation}

In conclusion, 
          \begin{equation}
          		\begin{array}[c]{lll}
                |B(0,p,p)|\le &\displaystyle C\left(\|v\|_{\mathcal{V}}+\|y_0\|_{L^2(\Omega)}\right)B(0,p,p)^{1/2}+M\beta\int_{Q}\varrho^{-2}|L_0^*(p)|^2\,dxdt+\left(M^2+\frac{M}{\beta}\right)S^2\int_{Q}\varrho_0^{-2}|p|^2\\
                &\le \displaystyle C\left(\|v\|_{\mathcal{V}}+\|y_0\|_{L^2(\Omega)}\right)B(0,p,p)^{1/2}+
                \max\left\{ \beta M,\left(M^2+\frac{M}{\beta}\right)S^2\right\}B(0,p,p).
				\end{array}
		\end{equation}
        
 Remember that $\beta M<1$, and   from inequality $\frac{1}{M}\frac{S^2M^2}{1-M^2S^2}<\beta$  we get $(M^2+M/\beta)S^2<1$ and the term  $\max\left\{ \beta M,\left(M^2+\frac{M}{\beta}\right)S^2\right\}B(0,p,p) $ can be absorbed to the left hand side to get 
     
    		\begin{equation}
    				\begin{array}[c]{lll}
   \|f[v]\|_{\mathcal{F}}+\|y\|_{\mathcal{Y}} <C\left(\|y_0\|_{L^2(Q)}+\|v\|_{\mathcal{V}}\right).
    				\end{array}
    		\end{equation}
            
\end{proof}

  \begin{lema}\label{pp674hg}
 Let be $\Phi[v]$ the set of all followers $f\in \mathcal{F}$  that solve Problem \ref{4.5.6}. Then the set $\mathcal{G}=\left\{(v;f):v\in \mathcal{U},f\in \Phi[v] \right\} $ is weakly closed. Moreover the functional $P(v;f)$ given  by \eqref{5.6.7} is lower semicontinous. 
  \end{lema}
  
 \begin{proof}

 First we proof   the weak closedness of $\mathcal{G}$. Let $(v_n; f_n)$  a sequence in $\mathcal{G}$ that converges to $(v;f)$. Suppose that $(v;f)$ is not in $\mathcal{G}$ so it exists a pair $(v;\tilde{f})$ such that
 
 			 	\begin{equation}
   			\frac{1}{2}\int_{Q}\varrho^2|\tilde{y}|^2+\frac{1}{2}\int_{\gamma \times (0,T)}\varrho^2|\tilde{f}|^2 < \frac{1}{2}\int_{Q}\varrho^2|y|^2+\frac{1}{2}\int_{\gamma \times (0,T)}\varrho^2|f|^2,
   	\end{equation} 
 where $\tilde{y}$ solves 
 
 			\begin{equation}
         				\begin{array}[c]{lll}
            							\tilde{y}_t-\Delta\tilde{y}+F(\tilde{y})=0 & \hbox{in}&Q,\\
                                        \tilde{y}=\tilde{f}1_{\gamma} +v1_\sigma & \hbox{on}&\Sigma,\\
                                        \tilde{y}(0)=y_0 & \hbox{in}&\Omega.
            					\end{array}
            			\end{equation}
 Observe that it exist $\delta >0$ such that 
 
 			\begin{equation}
   			\frac{1}{2}\int_{Q}\varrho^2|\tilde{y}|^2+\frac{1}{2}\int_{\gamma \times (0,T)}\varrho_0^2|\tilde{f}|^2 +\delta < \frac{1}{2}\int_{Q}\varrho^2|y|^2+\frac{1}{2}\int_{\gamma \times (0,T)}\varrho_0^2|f|^2.
   	\end{equation}
 Since $(v_n; f_n)$ converges to $(v;f)$  for some natural number $N$ and $n\geq N$ 
 
   	\begin{equation}
   			\frac{1}{2}\int_{Q}\varrho^2|\tilde{y}|^2+\frac{1}{2}\int_{\gamma \times (0,T)}\varrho_0^2|\tilde{f}|^2 < \frac{1}{2}\int_{Q}\varrho^2|y_n|^2+\frac{1}{2}\int_{\gamma \times (0,T)}\varrho_0^2|f_n|^2.
  	\end{equation}
    Taking the value of $S$ along the sequence $(v_n,\tilde{f})$ then it is possible to see that  
    	\begin{equation}
   			\frac{1}{2}\int_{Q}\varrho^2|\tilde{y}_n|^2+\frac{1}{2}\int_{\gamma \times (0,T)}\varrho_0^2|\tilde{f}|^2 \leq \frac{1}{2}\int_{Q}\varrho^2|y_n|^2+\frac{1}{2}\int_{\gamma\times (0,T)}\varrho_0^2|f_n|^2,
   	\end{equation}
  that contradicts the fact that $(v_n;f_n)\in \mathcal{G}$.  Let $\{(v_n; f_n)\}\subset \mathcal{G}$ be a sequence such that  $\|f_n\|_{\mathcal{F}}\rightarrow \infty$. Then, by inequality $ \|f_n\|_{\mathcal{F}}+\|y\|_{\mathcal{Y}} <C\left(\|y_0\|_{L^2(Q)}+\|v_n\|_{\mathcal{V}}\right)$ given in Proposition \ref{0.zq} , is straightforward to see that $\|v_n\|_{\mathcal{V}}\longrightarrow \infty$.  Then $P$ is coercive.  
  
  Finally, we prove that $P$ attains an infimum in $\mathcal{G}$. Let $(v_n; f_n)\in \mathcal{G}$ be a minimizing sequence. By definition, the sequence $(v_n; f_n)$ is uniformly bounded in $\mathcal{V}\times \mathcal{F}$   and then has  sub-sequence $(v_{n_k},f_{n_k})$ which is weakly convergent to some $(\hat{v}; \hat{f}[v]) \in \mathcal{G}$. This is the candidate to be a minimum.  Because $P$ is weakly lower semicontinous, it follows that
   		\begin{equation}
   				P(\hat{v}; \hat{f}[v])\leq \liminf_{n\rightarrow \infty}P(v_{n_k}; f_{n_k})\leq \inf_{(v;f)\in \mathcal{G}} P(v;f).
   		\end{equation}
   In conclusion,  $(\hat{v}; \hat{f}[v])$ is the solution to \eqref{func76638rtgh}.     
   
   \end{proof}
   
\noindent The next proposition addresses  the main objective of this paper, namely  Theorem  \ref{1.23} and give explicit solutions for the leader-follower pair.

  \begin{proposition}\label{OptS.s}
Let $T>0$ and $y_0\in L^2(\Omega)$ be given.  Then it exists a pair $(\hat{v};\hat{f}[\hat{v}])\in \mathcal{V}\times \mathcal{F}$ that is a solution to \eqref{5.6.7} and the associated state $\hat{y}$ solves the initial value problem
 
   \begin{equation} \label{a.9-1}
\left\{
\begin{array}[c]{lll} \displaystyle
\hat{y}_{t}-\Delta \hat{y} + F(\hat{y})=0 & \hbox{in} &Q,
\\  \displaystyle
\hat{y} = \hat{f}[\hat{v}]1_{\gamma}+ \hat{v}1_\sigma & \hbox{on}  & \Sigma,
\\  \displaystyle
\hat{y}(0) = y_{0}, \hat{y}(T)=0 & \hbox{in}&\Omega.
\end{array}
\right.
   \end{equation}
Moreover, it is possible to characterize 
   \begin{equation}\label{abx-2}
\hat{f}[\hat{v}] = \bigl.\varrho_0^{-2} \partial_\nu\hat{p}1_{\gamma }, \quad \hat{y} = \varrho^{-2} L_{F'(\hat{y})}^*( \hat{p}),
   \end{equation}
   where the function $\hat{p}\in \mathcal{P}$ solves the bilinear equation 
   
   				\begin{equation}
				\int_Q\varrho^{-2}L^*_{F'(\hat{y})}(\hat{p})L^*_{F_0(\hat{y})}(q)+\int_{\gamma \times(0,T)}\varrho_0^{-2}\partial_{\eta} \hat{p} \partial_{\eta} q\, d\Sigma\\
                =\int_{\sigma \times (0,T)} v \partial_\eta q\,d\Sigma+\int_\Omega y_0(x)q(0)dx,  \hspace{0.3cm} \forall q\in \mathcal{P}.
\end{equation}
  The leader control is expressed like
   
   \begin{equation}\label{ac.x-3}
\hat{v} = -\varrho^{-2}_*\partial_\eta (\hat{\zeta}+\hat{\phi})1_\sigma,
   \end{equation}
   where function $\hat{\zeta}$ is given by the initial value problem
   
   \begin{equation}\label{7.x-4}
\left\{
\begin{array}
[c]{lll}  \displaystyle
-\hat{\zeta}_{t} \!-\! \Delta \hat{\zeta} \!+\! F'(\hat{y})\hat{\zeta} 
\!=\! \alpha(\hat{y} \!-\! y_d) 1_{\Omega_d} \!-\! F'(\hat{y})\hat{\phi} \!-\! \varrho^{-2} F''(\hat{y}) \hat{p} L^*_0\hat{\phi}
& \hbox{in} &Q,
\\  \displaystyle
\hat{\zeta} = 0 & \hbox{on}  & \Sigma,
\\  \displaystyle
\hat{\zeta}(T) = 0 & \hbox{in}&\Omega,
\end{array}
\right.
   \end{equation}
 and the function $\hat{\phi}\in \mathcal{P}$ fulfills the equation
  
  \begin{equation}\label{au.j5}
 \displaystyle\int_{\gamma\times(0,T)}\varrho_0^{-2}\partial_\eta q \partial_\eta(\hat{\zeta}+\hat{\phi})\,d\Sigma+\int_Q\varrho^{-2}L_0^*(\hat{\phi}) L^*_{F'(y)}(q)\,dxdt =0, \hspace{0.5 cm}
\forall q \in \mathcal{P}.
   \end{equation}
\end{proposition}

\begin{proof}

\textbf{1.}The central idea is to formulate the optimization problem \eqref{5.6.7} as a problem of optimization with constraints. Define $\mathcal{X}=\mathcal{Y}\times\mathcal{F}\times\mathcal{V}\times\mathcal{P}$ and the operator    $M:\mathcal{X}\longrightarrow  L^2(Q)\times \mathcal{Y}\times\mathcal{F}$  given by 
   
   \begin{equation}
   		M(y,f,v,p)=\left(y+H_0(F(y))-G(f1_{\gamma}+v1_{\sigma}),y-\varrho^{-2}L^*_{F'(y)}(p),f+\varrho_0^{-2}\partial_{\eta}p1_{\gamma}\right),
   \end{equation}
where the operators $H_0$ and $G$ are given by \eqref{.fg3}. Define the functional

\begin{equation}\label{56.7}
			P_0(y,f,v,p)=\frac{\alpha}{2}\int_{Q_d}|y-y_d|^2+\frac{1}{2}\int_{\sigma\times(0,T)}\varrho_*^2|v|^2d\Sigma.
\end{equation}
Define the optimization problem given by

\begin{equation}\label{e.43}
\left\{
\begin{array}{l} \dis
\text{Minimize } P_0(y,f,v,p),
\\ \noalign{\smallskip} \dis
\text{Subject to } (y,f,v,p) \in \mathcal{X}, \
\\ \noalign{\smallskip} \dis \phantom{\text{Subject to } }
M(y,f,v,p) = (0,0,0).
\end{array}
\right.
   \end{equation}

\noindent By Lemma \eqref{pp674hg},  the above optimization problem has a solution $(\hat{y},\hat{f},\hat{v},\hat{p})\in \mathcal{X}$. The operator $M$ is linear and of class $C^1$  in the space $\mathcal{X}$ . Then given an arbitrary direction $(z,g,h,q)\in \mathcal{X}$  the derivative is given by
		\begin{equation}
				M'(y,f,v,p)(z,g,h,q)= \Big( z+H_0(F'(y)z)-G(g1_{\gamma}+h1_{\sigma}), z-\varrho^{-2}L_{F'(y)}(q)-F''(y)zp, g+\varrho_0^{-2}(\partial_{\eta}q) 1_{\gamma}\Big).
		\end{equation}
  Because $M$ is a closed operator, $M'(\hat{y},\hat{f},\hat{v},\hat{p})$ is closed, then $\text{ker} M'(\hat{y},\hat{f},\hat{v},\hat{p})^\bot=\text{Rank} M'(\hat{y},\hat{f},\hat{v},\hat{p})^*$. Define the operator $\mathcal{N}:\mathcal{P}\longrightarrow \mathcal{F},p\mapsto \varrho_0^{-2}\partial_\eta p$. By definition of the norm in $\mathcal{P}$ 
  $$\|\varrho_0^{-2}\partial_\eta p\|_{L^2(\Sigma)}\leq C\|p\|_{\mathcal{P}}$$ 
  so the operator $\mathcal{N}$ is continuous and then it  exists the dual operator $\mathcal{N}^*:\mathcal{F}\longrightarrow \mathcal{P}$. With this arguments one can get    
        		\begin{equation}\label{d.fg1}
        					 M'(\hat{y},\hat{f},\hat{v},\hat{p})^*(\psi,\phi,\varphi)	
                            \displaystyle=  \Big( \psi+\phi+F'(\hat{y})H_0^*(\psi)-\varrho^{-2}F''(\hat{y})\hat{p}\phi, \varphi-G^*(\psi)1_\gamma,-G^*(\psi)1_{\sigma},\mathcal{N}^*(\varphi1_{\gamma}) -L_{F'(\hat{y})}(\varrho^{-2}\phi)\Big).
        		\end{equation}

 \noindent Take the Gateaux derivate
		\begin{equation}
			\mathcal{D}P(\hat{y},\hat{f},\hat{v},\hat{p})(z,g,h,q) =\int_{Q_d}\alpha(\hat{y}-y_d)z\,dxdt +\int_{\sigma\times(0,T)}\varrho_*^2\hat{v}h\,d\Sigma,
		\end{equation}
Invoking Proposition \ref{dv.yu},  for $(\hat{z},\hat{g},\hat{h},\hat{q})\in  \hbox{Rank} M'(\hat{y},\hat{f},\hat{v},\hat{p})^*$, equation \ref{1.3gaw} takes the form  
			
            \begin{equation}\label{6.33}
  				\lambda\left(\alpha(\hat{y}-y_d)1_{Q_d},0,\varrho^2_*\hat{v},0\right)+ (\hat{z},\hat{g},\hat{h},\hat{q})=0.
  		\end{equation}
 Then, it exist   $(\psi,\phi,\varphi)\in L^2(Q)\times \mathcal{Y}\times\mathcal{F}$ such that 
    
 			\begin{equation}
 						\begin{array}[c]{lll}
 								\hat{z}&=\psi+\phi+F'(\hat{y})H^*(\psi)-\varrho^{-2}F''(\hat{y})\hat{p}\phi,\\
                                \hat{g}& = \varphi-G^*(\psi)1_\gamma,\\
                                \hat{h}&= -G^*(\psi)1_{\sigma},\\
                                \hat{q} & = \mathcal{N}^*(\varphi 1_{\gamma}) -L_{F'(\hat{y})}(\varrho^{-2}\phi). \\
 						\end{array}
 			\end{equation}
    Replace this  in equation \eqref{6.33} and normalizing by setting  $\lambda =1$, we have that
 
 		\begin{equation}\label{.tz1}
 				\begin{array}[c]{lll}
 							&\alpha(\hat{y}-y_d)1_{\Omega_d}=\psi+\phi+F'(\hat{y})H_0^*(\psi)-\varrho^{-2}F''(\hat{y})\hat{p}\phi,\\
                                &0= \varphi-G^*(\psi)1_\gamma,\\
                               &\varrho^2_*\hat{v}1_\sigma= -G^*(\psi)1_{\sigma},\\
                                 &0  = \mathcal{N}^*(\varphi1_{\gamma}) -L_{F'(\hat{y})}(\varrho^{-2}\phi).
 				\end{array}
 		\end{equation}
          Take the first equation in  (\ref{.tz1}). Define $\hat{\psi}=H_0^*(\psi)$,  $\hat{\phi}=-H_0^*(\phi)$  and $\hat{\zeta}=\hat{\psi}-\hat{\phi}$ . Then, using  the properties $L_0^*(\hat{\phi})=-\phi$  and $L_0^*(\hat{\psi})=\psi$, we get
         		
                \begin{equation}
                		\begin{array}[c]{lll}
                		\alpha(\hat{y}-y_d)1_{Q_d}&=
                         L_{F'(\hat{y})}^*(\hat{\zeta})+F'(\hat{y})\hat{\phi}+\varrho^{-2}F''(\hat{y})\hat{p}L_0^*(\hat{\phi}),\\
                        \end{array}
                \end{equation}
                
    It  is possible to write the above equation 
            
            \begin{equation}
            		 L_{F'(\hat{y})}^*\left(\hat{\zeta}\right)=\varrho^2(\hat{y}-y_d)1_{Q_d}-F'(\hat{y})\hat{\phi}-\varrho^{-2}F''(\hat{y})\hat{p}L_0^*(\hat{\phi}),
            \end{equation}
            
      with $\hat{\zeta}(T)=0$.   Recall equation (\ref{.tz1}) and use the fact that $\psi = L_0^*(\hat{\zeta}+\hat{\phi})$  and the equality $G^*( L_0^*(\hat{\zeta}+\hat{\phi}))=\partial_\eta (\hat{\zeta}+\hat{\phi})$  then
       	\begin{equation}
       					\hat{v}  =  \displaystyle-\varrho^{-2}_*\partial_\eta (\hat{\zeta}+\hat{\phi})1_\sigma
       	\end{equation}
       From the fourth equation of \eqref{.tz1} and for $q\in \mathcal{P}$ then 
        
        			\begin{equation}
        					\begin{array}[c]{lll}
        								\displaystyle\int_Q(\mathcal{N}^*(\varphi1_{\gamma}) -L_{F'(\hat{y})}(\varrho^{-2}\phi))q\,dxdt
                                        =0.
        					\end{array}
        			\end{equation}
                     But $\hat{\psi}=\hat{\zeta}+\hat{\phi}$ and from the second relation in \eqref{.tz1} is possible to get
                    
                    \begin{equation}
                    		\int_{\gamma\times(0,T)}\varrho_0^{-2}\partial_\eta q \partial_\eta(\hat{\zeta}+\hat{\phi})\,d\Sigma+\int_Q\varrho^{-2}L_0^*(\hat{\phi}) L^*_{F'(\hat{y})}(q)\,dxdt =0.
                    \end{equation}
        This concludes the proof.
\end{proof}

\subsection{Nonlinear gradient case}
In this section, we study the variant of the main problem where a non linear potential of the form $F(y,\nabla y)$ affects the system. Consider $F(x,w)$ of type $C^2(\mathbb{R}^{n+1})$ Lipschitz with $x\in \mathbb{R}$ and $w\in \mathbb{R}^n$  and all hypothesis established at the beginning. Define
        \begin{equation}\label{syst5gg}
				\left\{\begin{array}[c]{lll}
						\displaystyle y_t-\Delta y +F(y,\nabla y)=0  	& \hbox{in} &Q, \\
                        y = f1_{\gamma}+v1_{\sigma} & \hbox{on} & \Sigma,\\
                        y(0) = y_0 & \hbox{in} & \Omega.
 				\end{array}
                \right.
		\end{equation}
Consider the linear equation given by 
 \begin{equation}\label{syst5gg56}
				\left\{\begin{array}[c]{lll}
						\displaystyle y_t-\Delta y +
                        A\cdot \nabla y+ay=0 	& \hbox{in} &Q, \\
                        y = f1_{\gamma}+v1_{\sigma}  & \hbox{on} & \Sigma,\\
                        y(0) = y_0 & \hbox{in} & \Omega,
 				\end{array}
                \right.
		\end{equation}
with $A\in L^\infty(Q)^n$ and $a\in L^\infty(Q)$. Null controlability for the linear system \eqref{syst5gg56} comes in a similar way as Proposition 2.1. The formula for the follower control is similar to \eqref{8.xy} where $p\in \mathcal{P}$  is the solution to equation
\begin{equation}\label{555}
				B(a,A;p,q) = \int_{\sigma\times(0,T)}v \partial_{\eta}q\, d\Sigma+\langle y_0,q(0) \rangle_{L^2(\Omega)},\hspace{0.2cm} \forall q\in \mathcal{P}, 
		\end{equation}
where the bilinear $B(a,A;)$, as done in \eqref{8.gf},  is defined by the operator
\begin{equation}
    L^*_{a,A}(p)=-p_t-\Delta p+ap-\text{div}(Ap).
\end{equation}
and we can obtain bounds similar to \eqref{3,.,d}. 

Now we solve the non linear problem. Observe that is possible to linearize $F$ in the form
\begin{equation}
    F(y,\nabla y)=F_0(y)y+F_1(y)\nabla y
\end{equation}
where\footnote{Denote $F_x'$ and $F_w'$ as de derivates in each variable.} 
\begin{equation}
    F_0(y)=\int_0^1F'_x(s y,s\nabla y)ds,\hspace{1cm} F_1(y)=\int_0^1F'_w(s y,s\nabla y)ds.
\end{equation}
Let $z\in L^2(Q)$.  Then $F_0(z), F_1(z)\in L^\infty(Q)$.  By the results of Section 15, pg. 78 in \cite{LionsMagenes} the solution to the linearized equation
\begin{equation}\label{sygg56}
				\left\{\begin{array}[c]{lll}
						\displaystyle \partial_t y_z-\Delta y +
                        F_1(z)\cdot \nabla y_z+F_0(z)y_z=0 	& \hbox{in} &Q, \\
                        y_z = f_z[v]1_{\gamma}+v1_{\sigma}  & \hbox{on} & \Sigma,\\
                        y_z(0) = y_0 & \hbox{in} & \Omega,
 				\end{array}
                \right.
		\end{equation}
 satisfies $y_z\in H^{1/2,1/4}(Q)$. Using Schauder's fixed point theorem in the same way as in Proposition \ref{5.e.f}, we prove the existence of a null control for \eqref{syst5gg}.  Following the ideas of Step 3 in Proposition \ref{5.e.f} we conclude that the solution to \eqref{syst5gg} is $y=\rho_0^{-2}L^*_{F'_x(\hat{y},\nabla \hat{y}), F'_w(\hat{y},\nabla \hat{y})}(p)$, where $p\in \mathcal{P}$ satisfies
\begin{equation}\label{biy576cd}
				\int_Q\varrho^{-2}L^*_{F'_x(\hat{y},\nabla \hat{y}), F'_w(\hat{y},\nabla \hat{y})}(\hat{p})L^*_{F_0(\hat{y},\nabla \hat{y}), F_1(\hat{y},\nabla \hat{y})}(q )+\int_{\gamma \times(0,T)}\!\!\!\!\!\!\varrho_0^{-2}\partial_{\eta} \hat{p}\partial_{\eta} q\, d\Sigma\\
                =\int_{\sigma \times (0,T)} \!\!\!\!\!\!\!\!v\partial_{\eta}q\,d\Sigma+\int_\Omega y_0q_0, \hspace{0.3cm} \forall q \in \mathcal{P} .
		\end{equation}
To compute the leader control, follow the arguments in the proof of Proposition \ref{OptS.s} and we conclude with the next theorem.
\begin{theorem}
    Let $T>0$ and $y_0\in L^2(\Omega)$ be given. Then it exists a pair of controls  $(\hat{v},\hat{f}[\hat{v}])\in \mathcal{V}\times \mathcal{F}$  and the associated state $\hat{y}\in \mathcal{Y}$ solves \eqref{syst5gg}. We can characterize 
    \begin{equation}
       \hat{f}[\hat{v}] = \bigl.\varrho_0^{-2} \partial_\nu\hat{p}1_{\gamma }, \hspace{0.3cm} \hat{y}=\varrho^{-2}L^*_{F_x'(\hat{y}),F_w'(\nabla \hat{y})}(\hat{p}),
    \end{equation}
    where $\hat{p}\in \mathcal{P}$ solves \eqref{biy576cd}. The leader control has the form $\hat{v} = -\varrho^{-2}_*\partial_\eta (\hat{\zeta}+\hat{\phi})1_\sigma$  where functions $\hat{\zeta}$ , $\hat{\phi}\in \mathcal{P}$  solve the coupled system\footnote{ See that $DF$ denotes the derivate in $x,w$ and $F''_x$ and $F''_w$ second derivates. }
   
    \begin{equation}
\left\{
\begin{array}
[c]{lll}  \displaystyle
-L_{F_x'(\hat{y}),F_w'(\nabla \hat{y})}(\hat{\zeta}) 
\!=\! \alpha(\hat{y} \!-\! y_d) 1_{\Omega_d} \!-\!F'_x(\hat{y})\hat{\phi} \! - \hat{p}DF'_{x}(\hat{y})\cdot(I,\nabla)(\varrho_0^{-2}L_0^*(\hat{\phi}))
& \hbox{in} &Q,
\\ \displaystyle +\text{div}\left( F'_w(\hat{y},\nabla \hat{y})\hat{\phi}+ \hat{p}DF'_w(\nabla \hat{y})\cdot(I,\nabla)(\varrho_0^{-2}L_0^*(\hat{\phi}))\right)
\\  \displaystyle
\hat{\zeta} = 0 & \hbox{on}  & \Sigma,
\\  \displaystyle
\hat{\zeta}(T) = 0 & \hbox{in}&\Omega,
\end{array}
\right.
   \end{equation}
   
  \begin{equation}\label{au.j5g}
 \displaystyle\int_{\gamma\times(0,T)}\!\!\!\!\!\!\varrho_0^{-2}\partial_\eta q \partial_\eta(\hat{\zeta}+\hat{\phi})\,d\Sigma+\int_Q\varrho^{-2}L_0^*(\hat{\phi}) L^*_{F'_x(\hat{y})),F'_w(\nabla \hat{y})}(q)\,dxdt =0, \hspace{0.5 cm}
\forall q \in \mathcal{P}.
   \end{equation}

   \end{theorem}

\section{Inner leader control.}
In this section, the proof of Theorem \ref{1.23T} is given. 

   \begin{proposition}\label{0.2zq}
 Let $T>0$, $v\in \tilde{\mathcal{V}}$ and $y_0\in L^2(\Omega)$ .  It exists a follower control $f[v]\in \mathcal{F}$   such that  $y\in \mathcal{Y}$  is a solution to
  
   		\begin{equation}\label{2.5vc}
				\left\{\begin{array}[c]{lll}
						\displaystyle y_t-\Delta y +ay=v1_{\omega} 	& \hbox{in} &Q, \\
                        y = f[v]1_{\gamma}  & \hbox{on} &\Sigma,\\
                        y(0) = y_0, y(T)=0 & \hbox{in} & \Omega.
 				\end{array}
                \right.
		\end{equation}
        with   
 
 \begin{equation}\label{8.xc}
 f[v]=\varrho_0^{-2}\partial_\eta p1_{\gamma}; \hspace{0.5cm} y=\varrho^{-2}L_a^*(p),
 \end{equation}

 where $p\in \mathcal{P}$ solves
 		\begin{equation}\label{v.v.}
		B(a;p,q) = \int_{\omega\times (0,T)}vq+\langle y_0,q(0) \rangle_{L^2(\Omega)}, \hspace{0.3 cm} \forall q\in \mathcal{P}.
		\end{equation}
        Moreover

   \begin{equation}\label{7+<1}		\|f[v]\|_{\mathcal{F}}+\|y\|_{\mathcal{Y}}\leq C\left(\|v\|_{\tilde{\mathcal{V}}}+\|y_0\|_{L^2(\Omega)}\right).
    		\end{equation}.
\end{proposition}

\begin{proof}
			The proof  follows the  same strategy as done in Proposition \ref{0.zq}. We make several remarks. Optimization problem \ref{4.5.6}  equivalent to the fourth order system 
            \begin{equation}\label{2.5<x}
				\left\{\begin{array}[c]{lll}
						\displaystyle L_a(\varrho^{-2}L_a^*(p))=v1_{\omega} 	& \hbox{in} &Q, \\
                        \varrho^{-2}L_a^*(p) =-\varrho_0^{-2}\partial_\eta p1_{\gamma} & \hbox{on} &\Sigma,\\
                        (\varrho^{-2}L_a^*(p))(0) = y_0 & \hbox{in} & \Omega,
 				\end{array}
                \right.
		\end{equation}

where  $p\in \mathcal{P}$ solves the bilinear equation \eqref{v.v.}. Now we can estimate 
	\begin{equation}\label{5--a}		|l(q)|\leq\left(\int_{\omega\times(0,T)}\varrho_0^{2}|v|^2\right)^{1/2}\left(\int_{\omega\times(0,T)}\varrho_0^{-2}|q|^2\right)^{1/2}+\|y_0\|_{L^2(\Omega)}\|q(0)\|_{L^2(\Omega)}.
	\end{equation}
    
We proceed as in the previous case  to get inequality \eqref{7+<1}.

               \end{proof}

\subsection{Solution to the hierarchical control problem in the semi-linear case.}
 Given a function $z\in L^2(Q)$ it is possible to see that $F_0(z)\in L^\infty(Q)$. We consider the linearization of (\ref{syst123})   
    		\begin{equation}\label{syst57}
				\left\{\begin{array}[c]{lll}
						\displaystyle y_t-\Delta y +F_0(z)y=v1_{\omega} 	& \hbox{in} &Q, \\
                        y = f1_{\gamma}  & \hbox{on} &\Sigma,\\
                        y(0) = y_0 & \hbox{in} & \Omega.
 				\end{array}
                \right.
		\end{equation}

  \begin{proposition}\label{7.tg}
  Let $y_0\in L^2(\Omega)$, control $v\in \tilde{\mathcal{V}}$ and a positive time $T>0$. Then, it exists a follower control $f[v]\in \mathcal{F}$  and  $y\in \mathcal{Y}$ that solves the initial value problem 
   
   				\begin{equation}\label{sy7}
				\left\{	\begin{array}[c]{lll}
						\displaystyle y_t-\Delta y +F(y)=v1_{\omega} 	& \hbox{in} &Q, \\
                        y = f[v]1_{\gamma}  & \hbox{on} &\Sigma,\\
                        y(0) = y_0, y(T)=0 & \hbox{in} & \Omega,
 				\end{array}\right.
		\end{equation}

   		\begin{equation}\label{con57}
   				f[v]=\varrho_0^{-2}\partial _{\eta}p1_{\gamma}; \hspace{1cm}  y =\varrho^{-2}L^*_{F'(y)}(p),
   		\end{equation}
        where $p\in \mathcal{P}$ solves
      
      		\begin{equation}
      				\int_Q\varrho^{-2}L^*_{F'(y)}(p)L^*_{F_0(y)}(q)+\int_{\gamma \times(0,T)}\varrho_0^{-2}\partial_\eta p\partial_ \eta q\,d\Sigma=\int_{\omega\times (0,T)} v q+ \int_{\Omega} y_0q(0)\,dx, \hspace{0.3cm} \forall q\in \mathcal{P}.
      		\end{equation}
     Moreover it is possible to get the estimate
 
 \begin{equation}\label{9--.}
 \|f[v]\|_{\mathcal{F}}+\|y\|_{\mathcal{Y}}\leq C\left(\|v\|_{\tilde{\mathcal{V}}}+\|y_0\|_{L^2(\Omega)}\right). 
 \end{equation}
 
  \end{proposition}

\begin{proof}
\textbf{1.} Given $v\in \tilde{\mathcal{V}}$, we take $z\in L^2(Q)$ so $F_0(z)\in L^\infty(Q)$. By Proposition \ref{0.2zq} it exists    $f_z[v]\in \mathcal{F}$ such that  $y_z(T)=0$  where $y_z$ solves the equation 

		\begin{equation}\label{7.--}
				\left\{\begin{array}[c]{lll}
						\displaystyle y_{z,t}-\Delta y_z +F_0(z)y_z=v1_{\omega} 	& \hbox{in} &Q, \\
                        y_z = f_z[v]1_{\gamma}  & \hbox{on} &\Sigma,\\
                        y_z(0) = y_0, y_z(T)=0 & \hbox{in} & \Omega.
 				\end{array}
                \right.
		\end{equation}

Following the same arguments done in Step 1 of Proposition \ref{5.e.f}, the map  $\Lambda:L^2(Q)\longrightarrow L^2(Q)$, $z \mapsto y_z$, where $y_z$ solves \eqref{7.--}. is compact. Then it exist a solution $y\in \mathcal{Y}$ to \eqref{sy7}.

\textbf{2. } The existence of solutions for  problem \eqref{4.5.6} follows  the same argument of minimizing sequences  done in Step 2 of   Proposition \ref{5.e.f}.   Invoke the operators $H_0$ and $G$ from Proposition \ref{5.e.f}. Define the map $M:\mathcal{Y}\times \mathcal{F}\longrightarrow L^2(Q)$ given by 
	\begin{equation}
			M(y,f)=y-H_0(v1_{\omega}-F(y))-G(f1_{\gamma})-\bar{y}.
	\end{equation}
and the optimization problem with constrains 
 
 		\begin{equation}\label{prob54t}
 				\left\{
                \begin{array}[c]{lll}
 							\displaystyle \inf \frac{1}{2}\int_Q\varrho^2|y|^2+\frac{1}{2}\int_{\gamma\times(0,T)}\varrho_0^2|f|^2\,d\Sigma,\\
                            M(y,f)=0,\\
                            (y,f)\in \mathcal{Y}\times \mathcal{F}.
 				\end{array}
                \right.
 		\end{equation}
  
  	The derivative and the adjoint operators are given by
    
    \begin{equation}
  				M'(y,f)(z,g)=z+H_0\left(F'(y)z\right)-G(g1_{\gamma}).
  		\end{equation}
        
and
\begin{equation}
			M'(y,f)^*(\psi)=\left(\psi+F'(y)H_0^*(\psi), - G^*(\psi)1_\gamma\right).
	\end{equation}

	By Proposition \ref{dv.yu} we see that for $\psi\in L^2(Q)$, equation  \eqref{1.3gaw} takes the form

 	\begin{equation}\label{ehgjy6784}
 					-\lambda(\varrho^2y,\varrho_0^2f)+\left(\psi+F'(y)H_0^*(\psi),-G^*(\psi)1_\gamma\right)=0.
 	\end{equation}
 Define $p:=H_0(\psi)$.  Invoking similar computations as done before, we can reach the fourth order system 

                \begin{equation}\label{pro764hfj}
                		\left\{	\begin{array}[c]{lll}
                				\displaystyle L_{F_0(y)}(\varrho^{-2}L^*_{F'(y)}(p))=v1_{\omega}&\hbox{in}&  Q,\\
                                \varrho^{-2}L^*_{F'(y)}(p)=-\varrho_0^{-2}\partial_{\eta }p1_\gamma &\hbox{on}&  \Sigma,\\
                                \varrho^{-2}L^*_{F'(y)}(p)(0)=y_0, \varrho^{-2}L^*_{F'(y)}(p)(T)=0 &\hbox{in}&  \Omega.
                		\end{array}\right.
                \end{equation}

	Proceeding in the same way as in \eqref{ine.75648} and defining the same constants $M, S$ and $\beta$ we can achieve

		\begin{equation}\label{ineqhg75648}
				\begin{array}[c]{lll}
				B(0;p,p)&=\displaystyle \int_\Omega y_0(x)p(0)dx+\int_{\omega\times(0,T)}vp,dxdt\\
                &\displaystyle - \int_Q \varrho^{-2}\left(F_0(y)L^*_0(p)p+F'(y)L_0^*(p)p+F'(y)F_0(y)|p|^2\right)\,dxdt \\
                &\displaystyle \leq \left( \int_{\omega\times(0,T)}\varrho_0^{2}|v|^2 \right)^{\frac{1}{2}}\left( \int_{\omega\times(0,T)}\varrho_0^{-2}|p|^2 \right)^{\frac{1}{2}}+\left( \int_{\Omega}\varrho^{2}|y_0|^2 \right)^{\frac{1}{2}}\left( \int_{\Omega}\varrho^{-2}|p(0)|^2 \right)^{\frac{1}{2}}\\
                &\displaystyle - \int_Q \varrho^{-2}\left(F_0(y)L^*_0(p)p+F'(y)L_0^*(p)p+F'(y)F_0(y)|p|^2\right)\,dxdt \\
                & \displaystyle \leq \|y_0\|_{\mathcal{Y}}\max_{t\in[0,T/2]}\|p(t)\|_{L^2(\Omega)} +\|v\|_{\mathcal{U}}\|p\|_{\mathcal{P}}+M^2\int_{Q}\varrho^{-2}|p|^2 \,dxdt\\
                &\displaystyle + 2M\int_Q\varrho^{-2}|L_0^*(p)||p|^2\,dxdt.  \\
           \end{array}
           \end{equation}

    Then we get \eqref{9--.}.
        \end{proof}

\begin{lema}\label{v.gf}
  		 Let be $\Phi[v]$ the set of all followers $f\in \mathcal{F}$  that solve problem  \eqref{2..x} then the set $$\mathcal{G}=\left\{(v;f):v\in \tilde{\mathcal{V}},f\in \Phi[v] \right\} $$  is weakly closed in $\tilde{\mathcal{V}}\times \mathcal{F}$. Moreover the functional $P:\mathcal{G}\longrightarrow \mathbb{R}, (v;f)\mapsto P(v;f)$    is coercive and weakly lower semicontinous.  
  	\end{lema}
    
     We finish with the next proposition.
 \begin{proposition}\label{inner45}
 Let $y_0\in L^2(\Omega)$. It exists a pair $(\hat{v};\hat{f}[\hat{v}])\in \mathcal{G}$  such that the associated  $\hat{y}\in \mathcal{Y}$ solves 
 
 \begin{equation} \label{ab.51}
\left\{
\begin{array}[c]{l} \dis
\hat{y}_{t}-\Delta \hat{y} + F(\hat{y})=   \hat{v}1_\omega   \hbox{ in }  Q,
\\ \dis
\hat{y} = \hat{f}[\hat{v}]1_{\gamma}  \hbox{ on }    \Sigma,
\\ \dis
\hat{y}(\cdot,0) = y_{0}, y(0)=0   \hbox{ in } \Omega,
\end{array}
\right.
   \end{equation}
   \begin{equation}\label{an.2,}
\hat{f}[\hat{v}] = \bigl.\varrho_0^{-2} \partial_{\eta}\hat{p} \bigr|_{\gamma}, \quad \hat{y} = \varrho^{-2} L_{F'(\hat{y})}^* \hat{p}, 
   \end{equation}
\begin{equation}\label{26a.}
\int_Q \!\!\varrho^{-2} L^*_{F'(\hat{y})} \hat{p} \, L^*_{F_0(\hat{y})} p' \!+\! \int_{\gamma\times(0,T)} \varrho_0^{-2} \partial_{\eta} \hat{p}\partial_{\eta}p'\, d\Sigma 
= \int_{\omega\times (0,T)}\!\! \hat{v}   p' + \int_\Omega y_0(x) p'(x,0) \,dx \quad
\forall p' \in \mathcal{P}.
   \end{equation}
   \\
   Define $\hat{\gamma}$ the solution  to 
   
   \begin{equation}\label{aux-4.}
\left\{
\begin{array}
[c]{lll} \dis
-\hat{\gamma}_{t} \!-\! \Delta \hat{\gamma} \!+\! F'(\hat{y})\hat{\gamma} 
\!=\! \alpha(\hat{y} \!-\! y_d) 1_{\Omega_d} \!+\! F'(\hat{y})\hat{\phi} \!+\! \varrho^{-2} F''(\hat{y}) \hat{p} L^*_0\hat{\phi}
 &\hbox{ in }  Q,
\\ \dis
\hat{\gamma} = 0  &\hbox{ on }    \Sigma,
\\ \dis
\hat{\gamma}(\cdot,T) = 0  &\hbox{ in } \Omega,
\end{array}
\right.
   \end{equation}
   
   with  $\hat{\phi} \in \mathcal{P}$ the unique solution to
   
   \begin{equation}\label{aux-5.}
\int_{Q}\varrho^{-2}L_0^*(\hat{\phi})L_{F'(y)}^*(q)\,dxdt+\int_{\Sigma}\varrho^{-2}_0\partial_\eta  \hat{\phi} \partial_{\eta}q\, d\Sigma =- \int_{\Sigma}\varrho^{-2}_0\partial_\eta \hat{\gamma} \partial_\eta q \,d\Sigma
\quad \forall q \in \mathcal{P}. \ \
   \end{equation}
Then,  the leader control is characterized by
   \begin{equation}\label{aux.3}
\hat{v} = -\bigl.  \varrho_0^{-2} (  \hat{\gamma} +\hat{\phi} ) 1_{\omega } .
   \end{equation}

 \end{proposition}

\begin{proof}
Define $M:\mathcal{Y}\times \mathcal{U}\times\mathcal{F}\times \mathcal{P}\longrightarrow L^2(Q)\times \mathcal{Y}\times\mathcal{F}$   as 
	\begin{equation}
  			M(y,f,v,p)
  			=\big(y-H_0(v1_\omega-F(y))-G(f1_\gamma)-\bar{y},
  			y-\varrho^{-2}L^*_{F'(y)}(p),
  			f+\varrho_0^{-2}\partial_\eta p1_\gamma\big),
    \end{equation}
 where the operators $H_0$ and $G$  are the same as in\eqref{.fg3}. Define $P_0:M:\mathcal{Y}\times \mathcal{F}\times\tilde{\mathcal{V}}\times \mathcal{P}\longrightarrow \mathbb{R}$  
 
 \begin{equation}
 		P_0(y,f,v,p):=\frac{\alpha}{2}\int_{Q_d}|y-y_d|^2\,dxdt+\frac{1}{2}\int_{\omega\times(0,T)}\varrho_0^2|v|^2\,dxdt.
 \end{equation}

 Define the optimization problem
 
 	\begin{equation}\label{.,d}
 				\left\{\begin{array}[c]{lll}
 					\inf  P_0(y,f,v,p)  \\
                    M(y,f,v,p)=(0,0,0)\\
                    (y,f,v,p)\in \mathcal{Y}\times \mathcal{F}\times\tilde{\mathcal{V}}\times \mathcal{P}.
 			\end{array}\right.
 	\end{equation}

We have that
  	\begin{equation}
    	\begin{array}[c]{lll}
  			M'(\hat{y},\hat{f},\hat{v},\hat{p})(z,g,w,q)&=\Big(z-H_0(w1_\omega-F'(\hat{y})z)-G(g1_\gamma)-\bar{y},\\
            &z-\varrho^{-2}L^*_{F'(\hat{y})}(q)-\varrho^{-2}F''(\hat{y})z\hat{p},\\
            &g+\varrho_0^{-2}\partial_\eta q1_\gamma)\Big)
            \end{array}
  	\end{equation}
and the adjoint operator        
 	\begin{equation}
 			\begin{array}[c]{lll}
 					M'(\hat{y},\hat{f},\hat{v},\hat{p})^*(\psi,\phi,\varphi) &= \Big(\psi+F'(y)H_0^*(\psi)+\varphi-\varrho^{-2}F''(\hat{y})\hat{p}\varphi ,\\
                    &H_0^*(\psi)1_\omega,\\
                    &\varphi-G^*(\psi)1_\gamma,\\
                   & \mathcal{N}^* (\varphi 1_{\gamma})-L_{F'(\hat{y})}(\varrho^{-2}\phi)\Big).
 			\end{array}
 	\end{equation}

By Proposition \ref{dv.yu}  there exist functions$(\psi,\phi,\varphi)\in L^2(Q)\times \mathcal{Y}\times \mathcal{F}$  such that \eqref{1.3gaw}takes the form (as done several times before)

		\begin{equation}\label{eq56hgb}
				\left\{	\begin{array}[c]{lll}
						\alpha(\hat{y}-y_d)1_{\Omega_d}=\psi +F'(\hat{y})H_0^*(\psi)+\varphi-\varrho^{-2}F''(y)\varphi p \\
                         0=\varphi- G^*(\psi)1_\gamma\\
                      \varrho_0^2 \hat{v} =H^*_0(\psi)1_\omega \\
                        0=-L_{F'(y)}(\varrho^{-2}\phi)+\mathcal{N}^* (\varphi 1_{\gamma}).
				\end{array}\right.
		\end{equation}
        
Define the associated functions $\hat{\psi}=H_0^*(\psi)$, $\hat{\phi}=-H_0^*(\phi)$  and $\phi =L_0^*(H_0^*(\phi))=L_0^*(\hat{\phi})$. Taking $\hat{\gamma} =\hat{\psi}-\hat{\phi} $ and from  the first equality from equation \ref{eq56hgb} is possible to get

	\begin{equation}\label{gammeq}
			\begin{array}[c]{lll}
				\alpha(\hat{y}-y_d)1_{\Omega_d}=  L^*_{F'(\hat{y})}(\hat{\gamma})+F'(\hat{y})\hat{\phi}+\varrho^{-2}F''(\hat{y})L_0^*(\hat{\phi)}\hat{p}.
			\end{array}
	\end{equation}
    
  From \eqref{gammeq} we can obtain the equation
    
    		\begin{equation}
    				\begin{array}[c]{lll}
    						-\hat{\gamma}_t-\Delta\hat{\gamma}+F'(\hat{y})\hat{\gamma}=\alpha(\hat{y}-y_d)1_{\Omega_d}-F'(\hat{y})\hat{\phi}-\varrho^{-2}F''(\hat{y})L_0(\hat{\phi})\hat{p}\\
                            \hat{\gamma}(T)=0.
    				\end{array}
    		\end{equation}
           
    From \ref{eq56hgb} we can get
    	\begin{equation}\label{re76hgh}
    			\hat{v}=-\varrho^{-2}_0(\hat{\gamma}+\hat{\phi})1_{\omega}.
    	\end{equation}
    Taking fourth  equation  in \eqref{eq56hgb} and for $q\in \mathcal{P}$

    			\begin{equation}
        					\begin{array}[c]{lll}
        								\dis &\dis\int_Q(\mathcal{N}^*(\varphi1_{\gamma}) -L_{F'(y)}(\varrho^{-2}\phi))q\,dxdt\\
        								 \noalign{\smallskip}
  &=\dis\int_{\gamma\times(0,T)}\varrho_0^{-2}\partial_\eta q \partial_\eta\hat{\psi}\,d\Sigma+\int_Q\varrho^{-2}L_0^*(\hat{\phi}) L^*_{F'(y)}(q)\,dxdt\\
                                        &=0.
        					\end{array}
        			\end{equation}
                    
    Taking $\hat{\psi}=\hat{\phi}+\hat{\gamma}$
    
    		\begin{equation}
    				\int_{Q}\varrho^{-2}L_0^*(\hat{\phi})L_{F'(\hat{y})}^*(q)\,dxdt+\int_{\Sigma}\varrho^{-2}_0\partial_\eta \left(\hat{\phi}\right)\partial_{\eta}q\, d\Sigma = -\int_{\Sigma}\varrho^{-2}_0\partial_\eta\hat{\gamma}\partial_\eta(q)\,d\Sigma.
    		\end{equation}
            
   Then the theorem is proved.
\end{proof}

   \subsection{Fixed point iterations.}
We can find solutions to the hierarchical control problem via fixed point iterative methods in view of  Proposition \ref{inner45}. Denote by $l_n$ the linear functional  \eqref{v.v.} associated with $v_n$.
\begin{enumerate}
    \item Choose an initial $v_1\in \mathcal{V}$ and $y_1\in \mathcal{Y}$ ( with $f_1=0$) such that 
\begin{equation}
\begin{array}[c]{lll}
    \partial_ty_1-\Delta y_1+F(y_1)=v_1\mathbf{1}_{\omega}& \text{in}&Q\\
    y_1(0)=y_0 &\text{in} &\Omega.
    \end{array}
\end{equation} 
    \item  Let $n\in \mathbb{N}$. Define $p_{n}\in \mathcal{P}$  $$\int_Q \!\!\varrho^{-2} L^*_{F'(y_{n-1})} p_n \, L^*_{F_0(y_{n-1})} q \!+\! \int_{\gamma\times(0,T)} \varrho_0^{-2} \partial_{\eta} p_n\partial_{\eta}q\, d\Sigma 
= l_{n-1}(q), \hspace{0.5cm}\forall q\in \mathcal{P},$$ and  $f_{n}=f[v_{n-1}]\in \mathcal{F}$ and $y_{n}\in\mathcal{Y}$ via  \eqref{an.2,}. Moreover we have
 \begin{equation}\label{9--.n}
 \|p_n\|_{\mathcal{P}}+\|y_n\|_{\mathcal{Y}}\leq C\left(\|v_{n-1}\|_{\mathcal{V}}+\|y_0\|_{L^2(\Omega)}\right). 
 \end{equation}
    \item   Define $v_{n}\in \mathcal{V}$ via \eqref{aux.3}, where both $\phi_n,\gamma_n\in \mathcal{P}$ solve \eqref{aux-4.}.
\end{enumerate} 

 Coupled equations \eqref{aux-4.} allow to bound 
 \begin{equation}
 \|\phi_n\|_{\mathcal{P}} +\|\gamma_n\|_\mathcal{P}\leq C (\|p_n\|_{\mathcal{P}}+\|y_n\|_{\mathcal{Y}}).
 \end{equation}
   By linearity of the relation \eqref{aux.3}, we can find a constant $K<1$ (depending on $\alpha$) such that
 \begin{equation}
     \|v_{n}-v_{n-1}\|_{\mathcal{V}} \leq K\|v_{n-1}-v_{n-2}\|_{\mathcal{V}}.
 \end{equation}
 Then $v_n$ converges to some $v$ in $\mathcal{V}$ . The global Lipschitz continuity of $F$  implies strong convergence $F(y_n)\rightarrow F(y)$ and $F'(y_n)\rightarrow F'(y)$ . Since all auxiliary equations are linear, we may pass to the limit and conclude that  $(f,v,y,\phi,\gamma)$ solves the coupled equation of Proposition \ref{inner45}. These arguments can also be applied to Proposition 2.3.
 \section{Conclusions and further problems.}
The existence of controls for  \eqref{syst56ghj}  was  possible due the Lipschitz condition of the function $F$ that gives uniform bounds and  compactness of the map $\Lambda$.   Based on the results from \cite{fernandez1997null} we propose the hierarchical control problem with  super-linear behavior
        \begin{equation}
            \lim_{s\rightarrow\infty}\frac{F(s)}{s\log^\beta(1+|s|)}=0, \hspace{0.4cm} \beta >0.
        \end{equation}
  There are several results for  null controlability problem as  in \cite{deTeresaIns} and \cite{fernandez2000null} for the case $0<\beta< 3/2$.    Blowing-up phenomena for \eqref{syst56ghj} appear when $\beta>2$ and the remaining cases are still open. In \cite{fernandez2000null} it is explained that the presence of the term $1+1/T+(T^{1/2}+T)\|a\|+\|a\|^{2/3}$ in the observability inequality is the key idea for controlability  when  $0<\beta<3/2$.  For the hierarchical control problem proposed here, the presence of the fourth order system in Proposition \ref{0.zq} interfere with the sufficiency of this observability constant   to compensate the superlinear behavior of $F$.  It is an open problem the same problem in this setting.
  \\
  Also we propose the hierarchical control problem with Neumann boundary conditions given by 
\begin{equation}\label{nwm1}
				\left\{\begin{array}[c]{lll}
						\displaystyle y_t-\Delta y +F(y)=0  	& \hbox{in} &Q, \\
                        \partial_\eta y = f1_{\gamma}+v1_{\sigma} & \hbox{on} & \Sigma,\\
                        y(0) = y_0 & \hbox{in} & \Omega.
 				\end{array}
                \right.
		\end{equation}
Solving the above system requires new space of functions $p\in C^2(\bar{Q})$ where the normal derivates $\partial_ \eta p|_{\Sigma}=0$. This implies the need to compute new Carleman inequalities with boundary integrals of the form $\int_{\Sigma}\varrho^{-2}|p|^2d\sigma$  that, as the author concern, they are not deeply studied in the recent literature.  Therefore, problem is worth  exploring further and in general the mixed condition   $\partial_\eta y+\alpha y=f1_\gamma+v1_\sigma$ .

\appendix  
\section{Appendix: proof of Lemma \ref{lemmanm}}

\begin{proof}
Let  $q\in \mathcal{D}(\Sigma)$ a distribution.  Then 

		\begin{equation}\label{g8.x.}
				\int_{Q}L_0(G(q))p\,dxdt = \int_{\Sigma}p\partial_{\eta}G(q)\,d\Sigma-\int_{\Sigma}\partial_{\eta}pG(q)\,d\Sigma+\int_{Q}G(q)L_0^*(p)\,dxdt.
		\end{equation}    
   By the definition of an adjoint
   
    		\begin{equation}\label{.56z}
    				\int_{Q}G(q)L_0^*(p)\,dxdt=\int_{\Sigma}qG^*(L_0^*(p))\,d\Sigma. 
    		\end{equation}
      Then equation \eqref{g8.x.} together with \eqref{.56z} gets
      
      \begin{equation}\label{g8.q.}
				\int_{Q}L_0(G(q))p\,dx\,dt =\int_{\Sigma}p\partial_{\eta}G(q)\,d\Sigma -\int_{\Sigma}\partial_{\eta}pG(q)\,d\Sigma+\int_{\Sigma}qG^*(L_0^*(p))\,d\Sigma.
		\end{equation}
        
        Since $p|_\Sigma=0$ and from \eqref{.fg3} then (\ref{g8.q.}) becomes  
        
        \begin{equation}\label{6.v.9}
        		\int_{\Sigma}\partial_{\eta}pG(q)\,d\Sigma=\int_{\Sigma}qG^*(L_0^*(p))\,d\Sigma.
        \end{equation}
        Now $G(q)|_{\Sigma}=q|_{\Sigma}$ from \eqref{.fg3} and therefore for any $q$ in $L^2(\Sigma)$ the equality  (\ref{6.v.9}) takes the form
   		
   				\begin{equation*} 
       		\int_{\Sigma}q\partial_{\eta}p\,d\Sigma=\int_{\Sigma}qG^*(L_0^*(p))\,d\Sigma, 
        \end{equation*}
        
   	and then $G^*(L_0(p))=\partial_{\eta }p$.

    \end{proof}

\newpage

\bibliographystyle{unsrt}
\bibliography{biblio}

\end{document}